\documentclass[letterpaper, 10 pt, conference]{ieeeconf}  

\usepackage{graphicx}
\usepackage{epstopdf}
\usepackage{amsfonts}
\usepackage{amsmath,amssymb}
\usepackage{color}
\usepackage{tikz}
\usetikzlibrary{automata}
\usepackage{color}
\usepackage{verbatim}
\usepackage{algorithm}

\usepackage{fancyhdr}
\newtheorem{proposition}{Proposition}
\usepackage{varioref}
\let\labelindent\relax
\usepackage{enumitem}
\begin{document}
\title{\LARGE \bf A pedestrian hopping model and traffic light scheduling for pedestrian-vehicle
mixed-flow networks }
\author{Yi Zhang, Rong Su, Kaizhou Gao, Yicheng Zhang\\}

\maketitle

\begin{abstract}
This paper presents a pedestrian hopping model and a traffic signal scheduling strategy with consideration of both pedestrians and vehicles in the urban traffic system. Firstly, a novel mathematical model consisting of several logic constraints is proposed to describe the pedestrian flow in the urban traffic network and its dynamics are captured by the hopping rule, which depicts the changing capacity of each time interval from one waiting zone to another. Based on the hopping mechanism, the pedestrian traffic light scheduling problems are formulated by two different performance standards: pedestrian delay and pedestrian unhappiness. Then the mathematical technique and the meta-heuristic approach are both adopted to solve the scheduling problem: Mixed integer linear programming (MILP) formulation for pedestrian delay model and discrete harmony search algorithm (DHS) for both pedestrian delay model and unhappiness model. Secondly, a mathematical model about the vehicle traffic network, which captures drivers' psychological responses to the traffic light signals, is introduced. Thirdly, a traffic light scheduling strategy to minimize the trade-off of the delays between pedestrians and vehicles is proposed. Finally, we translate this traffic signal scheduling problem for both pedestrians and vehicles into a MILP problem which can be solved by several existing tools, e.g., GUROBI. Numerical simulation results are provided to illustrate the effectiveness of our real-time traffic light scheduling for pedestrian movement and the potential impact to the vehicle traffic flows by the pedestrian movement.
\end{abstract}


\section{Introduction} \label{sec:1}
In view of the increasing traffic congestion in urban road networks, tremendous efforts have been made to tackle this challenge via many different measures over the past several decades, and the urban traffic signal control is one of the most essential strategies among all these measures. The most existing studies on traffic signal control can be classified into two categories with four different types: fixed time strategies and traffic responsive strategies, isolated strategies and coordinated strategies. One of the well-known systematic mathematical frameworks for signal timing calculation was derived by Webster \cite{webster1958traffic}, which is only applicable to undersaturated conditions. The signal cycle was divided into separate stages and an empirical formula for the optimal cycle time to minimize overall junction delay was also given. SIGSET and SIGCAP proposed in \cite{allsop1971sigset} and \cite{allsop1976sigcap} are two typical examples of stage-based strategies, which all have stage-structure to be specified in the calculation. Phase-based strategies proposed in \cite{improta1984control} determines not only optimal splits and cycle time, which only consider in stage-based approach, but also the optimal staging. Group-based strategies deal with groups of streams without the need to maintain the stage-structure during optimization, which was firstly proposed by Gallivan and Heydecker \cite{gallivan1988optimising} and Heydecker and Dudgeon \cite{heydecker1987calculation}. Isolated traffic-responsive strategies proposed by Miller \cite{miller1963computer} use inductive loop detectors as real-time measurements to execute sophisticated vehicle-actuation logic. Fixed-time coordinated control strategies can be divided into two typical representatives: MAXBAND methods synchronize traffic signals so as to maximize the number of vehicles which can go through multi-intersections without stopping at any signal, see \cite{little1966synchronization} and \cite{gartner1981versatile}. TRANSYT firstly released in 1969 by Robertson \cite{robertson1900tansyt} is a computer model to optimize the linking and timing of traffic signals in a network. Since the settings of fixed-time strategies are based on historical data but not real-time data, traffic-responsive strategies make up the disadvantages once a real-time control system is installed. Corresponding Typical works are SCOOT, model-based optimization methods and store-and-forward based approaches. SCOOT is regarded as the traffic-responsive version of TRANYT and was first developed by Hunt and Robertson \cite{hunt1982scoot}, which run repeatedly in real time to investigate change of splits, offsets and cycle time at individual intersection. Many model-based traffic-responsive strategies have been released: OPAC\cite{gartner1983opac}, PRODYN\cite{henry1984prodyn}, CRONOS\cite{boillot1992optimal}, RHODES\cite{sen1997controlled}, which all solve a dynamic optimization problem in real time to obtain the optimal switching times. Store-and-forward model is a model simplification that deals with a mathematical description of the traffic flow process, and was first developed by Gazis and Potts \cite{gazis1963oversaturated}.\\
 \indent Both researchers and governments lay more emphasis on vehicle users, and the concentration on reducing delay for vehicular traffic leads to substantial delays to pedestrians, which has become increasing popular due to the recent trend of developing pedestrian friendly urban areas and the appearance of huge pedestrian volume in city centers. The installation and operation of a detection system (sensors, communications) has made it possible to further study pedestrian behavior and control pedestrian flow \cite{boudet2009pedestrian}. The unique characteristics that distinguish pedestrian traffic control have led to two different types of research focuses, i.e., safety-based studies and efficiency-based studies. In the safety-based studies, the existing literature studied the relationship between pedestrian fatalities and signal settings. Garder \cite{gaarder1989pedestrian} indicated that pedestrian safety analysis should be separated into two situations: (1) pedestrians and turning vehicles all pass the conflict area when the light is green, and (2) pedestrian non-compliance behavior occurs when crossing on red. A pedestrian accident prediction model for traffic signal was proposed, which analyzes the factors that influence pedestrian accident risks at signalized crosswalks in Poland \cite{tarko1995accident}. Moreover, traffic conflict technique (TCT), proposed by Perkins and Harris, determines the impact of pedestrian phase pattern on traffic safety \cite{perkins1968traffic}. While in the efficiency-based studies, many multi-objective analyses have been proposed, mainly based on either mathematical programming approaches or simulation-based methods. Optimal pedestrian signal timings were integrated into the corresponding optimization models either as one part of an objective function or as constraints, see \cite{allsop1992evolving} and \cite{lam1997integrated}. An optimization model, which selects pedestrian phase patterns between the normal two-way crossing (TWC) and the exclusive pedestrian phase (EPP) with the trade-off between safety and efficiency factors in an isolated intersection, is proposed in \cite{ma2014optimization} and \cite{ma2015optimization}. A simple hypothetical network with fixed-time noncoordinated signal cycles is analyzed for the effects of signal cycle timings on the delay caused by both vehicles and pedestrians \cite{ishaque2005multimodal}. Various multi-attribute weighting criteria are applied to traffic delays in order to examine cost trade-offs between pedestrians and vehicles \cite{ishaque2007trade}. The majority of the multi-objective analyses based on mathematical model only focus on one isolated intersection, rarely take the traffic network into account. Traffic signal control for pedestrian networks are seldom mentioned in the literature, which could be partially due to the uncontinuity of pedestrian flows in different intersections, in contrast to vehicle flows in an urban network. However, with the common adoption of GPS-enabled mobile devices for pedestrian applications, the establishment of pedestrian networks become possible by mobile users' volunteered data collections \cite{kasemsuppakorn2013pedestrian}. The current routes of navigation services are based on road networks and do not include off-road walking paths \cite{gaisbauer2008wayfinding}. Once pedestrian network databases are publicly available like road networks, the study on traffic signal control of pedestrian networks will be parallel with vehicle traffic control. \\
\indent Since the higher vehicle speed and the limited intersection length compared with the link length, it is reasonable to directly adopt cell transmission model (CTM) \cite{daganzo1994cell} on vehicle traffic network by neglecting the intersection length. However, as important as links in vehicle network, the intersections have the most vital role for pedestrian flows, which definitely cannot be neglected. Unlike the connected cells described in CTM of vehicle flows, a concept of "hopping" which describe the pedestrian flow from two isolated cells (waiting zones connected by the crosswalk) is firstly proposed. In order to cooperate with the traffic light assignments, our pedestrian hopping model guarantee the clearance of the pedestrian flow at the end of each GREEN sub-sequence. Then a network-based urban pedestrian traffic
light control problem is firstly formulated as a scheduling problem, aiming to reduce the total waiting time and the total unhappiness over a
time horizon. The pedestrian unhappiness is formulated as an exponential function in order to bring fairness to some pedestrians with few quantities but waiting for a long time, which frequently occurs in delay model. The dynamic of the pedestrian volume in each junction captures the characteristics of the conventional two-way crossing (north-south and east-west) of the crosswalk topology in each intersection. We assume that the pedestrian demands at each waiting zone are known in advance, and the diversion ratio, which decide their destination directions, and the capacity of each crosswalk are also known. These assumptions will be relaxed in our future research by using real-time data-driven model identification. Secondly, we convert the traffic light scheduling problem on pedestrian delay into a mixed integer linear program by using standard optimization tool, GUROBI \cite{optimization2012gurobi}. To overcome the huge computational complexity in mathematical technique, the discrete harmony search algorithm (DHS) is adopted to address the scheduling problems for both pedestrian delay and unhappiness, numerical experiments are presented to illustrate the effectiveness of our scheduling strategies, the fairness of our unhappiness model and the increased feasibility of obtaining real-time solution. Thirdly, a vehicle-based traffic light scheduling problem \cite{zhang2015urban}, aiming to reduce the total vehicle waiting time over a given time horizon, will also be introduced. Finally, we integrate our pedestrian model with the vehicle network model to formulate a linear optimization problem by developing a single objective function with different weight factors for each side, aiming to investigate the impact of pedestrians on vehicle delays when simultaneously providing pedestrians with convenient. Since both vehicle and pedestrian traffic light scheduling problems are captured as model-based optimization methods which do not consider explicitly splits, offsets, or cycles, the integration becomes much easier. Standard optimization tools such as GUROBI \cite{optimization2012gurobi} are also used to solve our integrated optimization problem. Numerical experiments are presented to illustrate the effectiveness of scheduling strategies on how the total vehicle delay is affected by pedestrians. \\
\indent This paper is organized as follow. The specific descriptions of a pedestrian hopping model and a formulation of traffic light scheduling on both pedestrian delay and pedestrian unhappiness are illustrated in Section \ref{sec:2}. The algorithms to solve the pedestrian optimization problems and the experimental results are also described in Section \ref{sec:2}. Then a vehicle traffic system model and the integration of both pedestrian delay and vehicle delay are presented in Section \ref{sec:3}, and the case study of integrated model is also described in Section \ref{sec:3}. Conclusion are drawn in Section \ref{sec:4}.

\begin{figure}[!ht]
	\centering
	\includegraphics[width=2.0in]{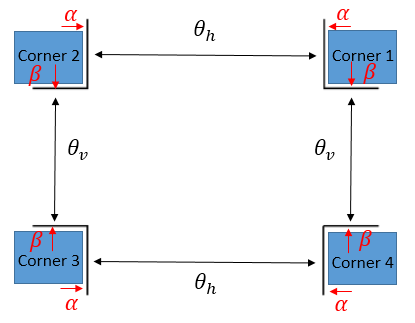}
	\caption{Pedestrian traffic light at junction $J$}
	\label{fig_sim2}
\end{figure}

\section{Formulation of A Pedestrian Traffic Light Scheduling Problem}\label{sec:2}

\subsection{A pedestrian hopping model} \label{subsecped}
An urban traffic network consists of a set of road links $L$ $\in \mathcal{L}$ and junctions $J$ $\in \mathcal{J}$. A simple pedestrian traffic network can be illustrated as a set of distributed subsystems (each junction) with four shaded areas, and these shaded areas are the assumptive waiting zones (corners shown in Fig. \ref{fig_sim2}) provided for pedestrians who want to cross the road, we use $W \in \mathcal{W}$ to denote the four waiting zones in each junction. Because of the uncontinuity of the relationship between the pedestrian demands at different junctions, which is different from a vehicle traffic network connected by incoming flow and outgoing flow from corresponding links, we assume that the pedestrian volume dynamics in one junction will not be influenced by those of other junctions. Let $\Omega_{p}^{J}$ be the set of pedestrian stages in junction $J$ (two stages: horizontal GREEN or vertical GREEN), $\mathcal{F}_{p}^{J}$ the set of pedestrian streams in junction $J$: $\mathcal{F}_{p}^{J} \in  \mathcal{W} \times \mathcal{W}$, and $h_{p}^{J}$ the association of each stage to relevant compatible streams, $h_{p}^{J}: \Omega_{p}^{J}\to 2^{\mathcal{F}_{p}^{J}}$. Fig. \ref{fig_sim2} depicts pedestrian traffic lights at junction $J$, where $\theta_{h}$ and $\theta_{v}$ denote the stages of pedestrian traffic lights at horizontal and vertical crosswalks respectively, and $\alpha$ and $\beta$ denote the diversion ratios at the corresponding corners. To simplify our presentation, we make the following assumptions:

	\begin{enumerate}[label=\textbf{\arabic*},start=1]
			\item Pedestrian phase patterns only consider the conventional two-way crossing, thus, exclusive pedestrian phases (e.g., diagonal crossing when vehicle traffic streams are stopped from all directions) are not considered.
			\item The pedestrian demand and diversion ratios ($\alpha$ and $\beta$) are known.
            \item The summation of diversion ratios $\alpha$ and $\beta$ at the same time instant are equal to 1 in one specific corner.		 
        	\item The capacity of a crosswalk during each time interval $\triangle$ is known.
            \item There are no pedestrians remaining on the crosswalk at the end of each green-time interval.
            \item Waiting-zone cells have infinite capacity.
		\end{enumerate}
Assumption 5 can be easily relaxed in our model, and we give this assumption here to better serve our case study in \ref{sec:3}. The detail of the relaxation is discussed in Section \ref{subsecintegration}.

\subsubsection{Parameters}
To better describe the pedestrian model, the parameters used in later part are summarized in TABLE \ref{label11}.

    \begin{table}[!ht]
		\caption{Notations for the pedestrian model}
       \label{label11}
		\centering
          \begin{tabular}{p{0.22\linewidth}p{0.7\linewidth}}
		  \hline
		  \textbf{Parameter} & \textbf{Definition}\\
		  \hline
            $i$ &  The index of waiting zones (corners) in the anticlockwise direction at a junction as shown in Fig. \ref{fig_sim2}.\\
            $j_{h}$ &  The index of adjacent corners that provide pedestrian \\
            ($j_{v}$) & flows in horizontal (vertical) directions.\\
            $o$ &  Pedestrian signal stages: horizontal direction or vertical direction, $o \in \Omega_{p}^{J}$.\\
            $\theta_{o}(k)$  &  The pedestrian traffic light in associated stage $o$. \\
			$f_{ij}(k)$  & The flow rate from corner $i$ to corner $j$ during time interval $k$, $i \in \mathcal{W}$, $j \in \mathcal{W}$.\\
           $P_{i}(k)$ & The number of pedestrians waiting to cross the street in corner $i$ at time instant $k$. \\
           $\hat{P}_{o}(k)$ &  The capacity of the crosswalk associated with stage $o$ during time interval $k$.\\
           $s_{i}(k)$ & The number of outgoing pedestrians of corner $i$ during time interval $k$. \\		
           $d_{i}(k)$ & The number of incoming pedestrians of corner $i$ during time interval $k$.\\
           $I_{i}(k)$ & The number of incoming pedestrians of corner $i$ without using crosswalk during time interval $k$,
            $I_{i}(k)$ is part of $d_{i}(k)$. \\
		   $\eta_{i}(k)$ & The pedestrian diversion ratio from corner $i$ to corner $j$ at time instant $k$, $\eta \in \mathcal{H} = \{\alpha, \beta\}$, where $\alpha$ is the \\ & horizontal diversion ratio from corner $i$ to corner $j_{h}$, $\beta$ is the vertical diversion ratio from corner $i$ to corner $j_{v}$, and $u_{p}^{J}: \Omega_{p}^{J} \to \mathcal{H}_{p}^{J}$ depicts the association of each stage to the pedestrian diversion ratio with corresponding direction, namely, if $o = h$, then $\eta = \alpha$, otherwise, if $o = v$, then $\eta = \beta$.\\
           $\gamma_{i}(k)$ & Pedestrians departure ratio at corner $i$ at time instant $k$, aiming at those who have already used crosswalk and reached destination.\\

\hline
		\end{tabular}
	\end{table}
 \subsubsection{Stage Constraints}
Due to the vehicle flow existing in the traffic network, the horizontal pedestrian traffic light cannot be green simultaneously with the vertical pedestrian signal.

    	\begin{subequations}\label{eq:stageped}
        \begin{align}
           &\underset{o \in \Omega_{p}^{J}}\sum\theta_{o}(k)=1 \label{eq:stageped1a}\\
           &(\forall w \in \Omega_{p}^{J}) \theta_{o}(k)=0 \nonumber \\
           & \Rightarrow (\forall (i,j) \in  h_{p}^{J}(o)) f_{ij}(k)=0 \label{eq:stageped1b}\\
           &(\forall o \in \Omega_{p}^{J})(\forall k \in \mathbb{N}) \theta_{o}(k)\in \{0, 1\} \label{eq:stageped1c}
        \end{align}
	    \end{subequations}
Where $\mathbb{N}$ denotes the set of natural numbers, $\theta_{o}(k) = 0$ and $\theta_{o}(k) = 1$ denote the RED and
GREEN traffic lights associated with stage $o$ respectively. Equation (\ref{eq:stageped1b}) illustrates that all associated flows with stage $o$ are zero if the stage traffic light is RED.

\subsubsection{Volume Dynamics}
    Let $i$ ($i \in \mathcal{W}$) be the index of corners in the anticlockwise direction at a junction as shown in Fig. \ref{fig_sim2}, namely, $i$ = \{1, 2, 3, 4\}, and
    $j_{h}$ ($j_{v}$) is the index of adjacent corners that provide pedestrian flows in horizontal (vertical) directions, and $j_{h}$ and $j_{v}$ follow the rule below:\\
    $j_{h} = \left\{\begin{array}{ll}
    i+1 & i \in \{1,3\}\\
    i-1 & i \in \{2,4\}
    \end{array}
    \right.$
    $j_{v} = \left\{\begin{array}{ll}
    j_{h}+2 & i \in \{1,2\}\\
    j_{h}-2 & i \in \{3,4\}
    \end{array}
    \right.$\\
The sources of pedestrians at corner $i$ are classified into three types: (1) pedestrians from its adjacent corner $j_{h}$, (2) from its adjacent corner $j_{v}$ and (3) freshmen from other places who have not crossed the road at this junction and have just entered the shaded corner $i$ waiting to cross. Similarly, the outputs of pedestrians at corner $i$ are also classified into three types: (1) pedestrians who want to go to corner $j_{h}$, (2) who want to go to corner $j_{v}$ and (3) those who have already crossed, reached their destination and do not want to cross anymore, and the third component of outputs is part of the first two components of sources which is specified by the departure ratio $\gamma$. The properties above can be captured by the following equations:

    	\begin{subequations}\label{eq:dyn1}
        \begin{align}
       (\forall k \in \mathbb{N})&P_{i}(k+1) = P_{i}(k) + \triangle(d_{i}(k)-s_{i}(k)) \label{eq:dyn1a}\\
        d_{i}(k)&= I_{i}(k) +  \underset{i \in \mathcal{C}:(j,i)\in \cup_{J \in \mathcal{J}}\mathcal{F}_{p}^{J}}\sum{f_{ji}(k)}  \qquad \label{eq:dyn1b}\\
        s_{i}(k)& =  \underset{i \in \mathcal{C}:(i,j)\in \cup_{J \in \mathcal{J}}\mathcal{F}_{p}^{J}}\sum{f_{ij}(k)}+ \nonumber  \\
        &\left \lfloor{ \gamma_{i}(k) \underset{i \in \mathcal{C}:(j,i)\in \cup_{J \in \mathcal{J}}\mathcal{F}_{p}^{J}}\sum{f_{ji}(k)} } \right \rfloor \label{eq:dyn1c}
        \end{align}
	    \end{subequations}

 \subsubsection{The Crossing Flow Constraints}
 The vehicle traffic network can be developed by adopting Cell Transmission Model, and the outgoing flow of the current link is usually captured by the current link volume, the remaining space of the downstream link and the capacity from current link to downstream link without considering the length of the intersection in many literatures \cite{lo1999novel}\cite{lin2001enhanced}\cite{han2014continuum}, which results from the high speed of the vehicles.
 Unlike roadways where vehicle flows connected by links, the topology structure of the pedestrian waiting zones around the intersection are four isolated areas connected by four crosswalks, which definitely cannot be neglected like macroscopic vehicle models since the relative lower speed of the pedestrians. Additionally, the essence of our adaptive control strategy is a real-time signal timing optimization algorithm, the traditional concepts of cycle time, splits and offsets, which are inherent in exiting signal optimization methods, will not appear in our adaptive algorithm. Instead, the optimal policy, a set of GREEN and RED sub-sequences with various lengths, is calculated for the entire prediction horizon in order to minimize our cost function, but only implemented the head section without degrading the performance of the optimization procedure.
 The optimized switching times of the GREEN-RED sequences are based on our pedestrian flow model, and we must guarantee that pedestrians on crosswalks should all be cleared at the end of each GREEN sub-sequence, and each GREEN sub-sequence is composed of several continuous GREEN time intervals.
 In view of the two problems illustrated above, the traditional cell transmission model for vehicle flows cannot be directly adopted in our pedestrian side, so the concept of the hopping, which satisfy both the topology structure of the pedestrian waiting zones and the clearance mechanism at the end of each GREEN sub-sequence, will be clearly discussed in the following content.\\
\indent The empirical model proposed by Virkler et al. \cite{virkler1998scramble} conducted 304 observations for two-way crossings to investigate the crossing time of pedestrian platoons, and they found that the crossing time for the front of the platoon is not affected by either the platoon size or the size of the opposing platoon while the crossing time for the rear of a platoon is affected by the size of the primary platoon but not significantly affected by the opposing platoon size. Moreover, the crossing time for pedestrian platoons, which is also recorded in HCM 2010 \cite{manual2010hcm2010}, was found after several regression analysis:

\begin{equation} \label{eq:crossingtime}
    T = \left\{
        \begin{array}{lll}
         I + \frac{L}{S_{p}}+ 0.27N_{ped}, & \text{if} \quad W \leq 3 m\\
         &\\
         I + \frac{L}{S_{p}}+0.81\frac{N_{ped}}{W}, & \text{if} \quad W > 3 m\\
        \end{array}
    \right.
    \end{equation}
where $T$ is the total crossing time ($s$), $I$ is the start-up time ($s$) and usually is 3.2s, $L$ is the length of crosswalk ($m$), $W$ is the width of crosswalk ($m$), $S_{p}$ is the average walking speed ($m/s$), and $N_{ped}$ is the size of the pedestrian platoon. The equation (\ref{eq:crossingtime}) indicate that the minimum time for $N_{ped}$ pedestrians successfully reaching the other curb under a certain length and width of the crosswalk. We can slightly change the form of the equation (\ref{eq:crossingtime}):

\begin{equation} \label{eq:capacityhop}
    N_{ped} = \left\{
        \begin{array}{lll}
         \left \lfloor{\frac{T - I - \frac{L}{S_{p}}}{0.27}}\right \rfloor, & \text{if} \quad W \leq 3 m\\
         &\\
         \left \lfloor{\frac{(T - I - \frac{L}{S_{p}})W}{0.81}} \right \rfloor, & \text{if} \quad W > 3 m\\
        \end{array}
    \right.
    \end{equation}
The equation (\ref{eq:capacityhop}) indicate that the maximum number of pedestrians $N_{ped}$ that can successfully finish crossing under a certain length and width of the crosswalk and a certain time interval $T$. Owing to the discrete characteristics of our model, we can set the sampling interval period $\triangle$ as our crossing time $T$, then corresponding capacity for pedestrian successfully finish crossing during this time interval can be obtained by equation (\ref{eq:capacityhop}) according to the state of the traffic light:\\
    \begin{enumerate}[label=\textbf{\arabic*},start=1]
    \item If $W \leq 3m$,
    \begin{equation} \label{eq:capacityfun1}
        \hat{P}_{o}(k) = \left\{
            \begin{array}{lll}
                &\left \lfloor{\frac{\triangle - I - \frac{L}{S_{p}}}{0.27}}\right \rfloor, \\
                & \qquad \quad \text{if} \quad (k = 1) \vee (\theta_{o}(k-1)=0)\\
                &\left \lfloor{\frac{\triangle}{0.27}}\right \rfloor, \\
                & \qquad \quad \text{if} \quad \theta_{o}(k-1)=1
            \end{array}
            \right.
    \end{equation}

    \item If $W > 3m$,
    \begin{equation} \label{eq:capacityfun2}
        \hat{P}_{o}(k) = \left\{
            \begin{array}{lll}
                &\left \lfloor{\frac{(\triangle - I - \frac{L}{S_{p}})W}{0.81}}\right \rfloor, \\
                & \qquad \quad \text{if} (k = 1) \vee (\theta_{o}(k-1)=0)\\
                & \left \lfloor{\frac{\triangle W}{0.81}}\right \rfloor, \\
                & \qquad \quad \text{if} \quad \theta_{o}(k-1)=1
            \end{array}
            \right.
            \end{equation}
\end{enumerate}

The equation (\ref{eq:capacityfun1}) and (\ref{eq:capacityfun2}) illustrate that the capacity for successfully crossing pedestrians in current GREEN time interval depends on the previous time interval: if the previous time interval is RED, then the capacity is a small value, otherwise, the capacity is a large value. The reason why capacity during one time interval $\triangle$ changes is the start-up time $I$ and the time spending on crosswalk $\frac{L}{S_{p}}$ appeared in equations (\ref{eq:crossingtime})-(\ref{eq:capacityfun2}). When the previous traffic signal state is RED, namely, $\theta_{o}(k-1)=0$, there are no pedestrians on the crosswalk at the start of current GREEN time interval ($\theta_{o}(k)=1$), however, when the previous traffic signal state is GREEN ($\theta_{o}(k-1)=1$), pedestrians are walking on the crosswalk at  the start of current GREEN time interval if there are sufficient demand on the waiting zone, accordingly, the capacity for the latter must be larger than the former. Therefore, the specific hopping flow can be described as:

\begin{equation} \label{eq:hopflow}
    f_{ij}(k)\triangle = \underset{\underset{\eta \in u_{p}^{J}(o)}{o \in \Omega_{p}^{J}: (i,j) \in  h_{p}^{J}(o)}}{\text{min}}(\hat{P}_{o}(k), P_{i}(k)\eta_{i}(k))
    \end{equation}
The equation above indicate that the outgoing hopping flow is determined by the capacity $\hat{P}_{o}(k)$ of current time interval and the demand in corner $i$ at time instant $k$. According to Assumption 6,  the waiting-zone cells can accommodate infinite volumes, so the remaining space of the waiting zone $j$ is not considered in our flow equation ($\ref{eq:hopflow}$). Because our signal scheduling is an on-line optimization strategy, the optimized signal setting can be obtained in advance, so the Flashing GREEN can be activated several seconds before the switching time from GREEN to RED. And we do not consider pedestrian non-compliance behavior in this paper, namely, pedestrians in waiting zones will all stop to cross once the traffic light turns Flashing GREEN. Therefore, pedestrians have already been cleared at the end of each GREEN sub-sequence, and directly switching to RED will not cause the safety issues.

\subsubsection{Objectives of the optimization model}
Based on the operational efficiency of the optimization model, we propose two different ways to measuare pedestrians' performance by determining the optimal signal timing for a signalized pedestrian network. The objectives of the optimization model are to minimize the total cost of delay and the total cost of unhappiness in entire network, respectively.
\paragraph{Linear Optimal Control for Pedestrian Delay}\label{sec:delaymodel}
Our cost function to minimize the total pedestrian delay time for the whole network within $H_{p}$ time intervals can be calculated as follows:
    \begin{equation}\label{eq:costped}
    \begin{aligned}
	\text{min} P_{D}&= \text{min}\sum_{J\in \mathcal{J}}\sum_{k=1}^{H_{p}}\sum_{i \in \mathcal{W}} [P_{i}^{J}(k)\\
    &- \underset{i \in \mathcal{C}:(i,j)\in \cup_{J \in \mathcal{J}}\mathcal{F}_{p}^{J}}\sum{f_{ij}(k)}]\triangle
    \end{aligned}
	\end{equation}
where $P_{i}^{J}(k)$ is the number of pedestrians of corner $i$ for intersection $J$ at time instant $k$, $\triangle$ is the sampling interval period.
\paragraph{Nonlinear Optimal Control for Pedestrian Unhappiness} \label{sec:unhappymodel}
Directly adopting delay time as a performance index may bring unfairness for some pedestrians. These pedestrians have the same destination and wait for a long time but still have no GREEN response, simply because the number of these people is so low that their sacrifices can maximize the global interest in entire network for a long time horizon. However, this is certainly unfair for these people, so a more sensible approach to take the urgent needs of this part of pedestrians into account is extremely necessary for the entire pedestrian network. To overcome such obstacle, an exponential function to express the extent of pedestrians' unhappiness is proposed, which is depicted by the duration of each RED sub-sequence and the number of RED sub-sequences over the whole prediction horizon.
    \begin{equation}\label{eq:costpedunhappy}
    \begin{aligned}
	\text{min} P_{U}&= \text{min}\sum_{J\in \mathcal{J}}\sum_{k=1}^{H_{p}}\sum_{o\in \Omega_{p}}\sum_{i\in \mathcal{W}} \bar{P}_{i,o}^{J}(k)\exp(\varphi_{o}^{J}(k)\triangle)
    \end{aligned}
	\end{equation}
where $H_{p}$ is the prediction horizon, $\triangle$ is the sampling interval period, and parameters $\varphi_{o}^{J}(k)$ and $\bar{P}_{i,o}^{J}(k)$ are clearly explained blow.\\
\indent Parameter $\varphi_{o}^{J}(k)$ can be obtained by introducing auxiliary functions $h_{o}^{J}(k)$, $f_{o}^{J}(k)$ and $q_{o}^{J}(k)$:

\begin{equation}\label{eq:h}
    h_{o}^{J}(k) = \left\{
        \begin{array}{lll}
            0, & \text{if} \quad k = 0 \\
            k [\theta_{o}^{J}(k)\oplus \theta_{o}^{J}(k+1)], & \text{if} \quad 1\leq k \leq H_{p}-1\\
            H_{p}, & \text{if} \quad k = H_{p}
        \end{array}
    \right.
\end{equation}

\begin{equation}\label{eq:f}
    f_{o}^{J}(k) = \left\{
        \begin{array}{lll}
            &max(\theta_{o}^{J}(k+1)-\theta_{w(J)}(k),0),  \\
            &\qquad \quad \text{if} \quad 1\leq k \leq H_{p}-1 \\
            &0, \qquad \text{if} \quad \theta_{o}^{J}(H_{p})=1\\
            &1, \qquad \text{if} \quad \theta_{o}^{J}(H_{p})=0
        \end{array}
    \right.
\end{equation}

\begin{equation}\label{eq:q}
    q_{o}^{J}(k) = \left\{
        \begin{array}{ll}
        &max(h_{o}^{J}(k)-h_{o}^{J}(k-1),0),\\&\qquad \qquad \qquad \qquad \qquad \text{if} \quad  h_{o}^{J}(k-1)\neq0\\
        &max(h_{o}^{J}(k)-h_{o}^{J}(k-1)-\sum_{i=0}^{k-1}q_{o}^{J}(i),0), \\
        & \qquad \qquad \qquad \qquad \qquad \text{if} \quad  h_{o}^{J}(k-1)=0
        \end{array}
    \right.
\end{equation}

Then the parameter $\bar{P}_{i,o}^{J}(k)$ associated with corner $i$ and stage $o$ at intersection $J$ in equation (\ref{eq:costpedunhappy}) can be obtained by introducing auxiliary function $\tilde{P}_{i,o}^{J}(k)$:
\begin{equation}\label{eq:P1}
    \tilde{P}_{i,o}^{J}(k) = \left\{
        \begin{array}{lll}
           &\underset{\eta \in u_{p}^{J}(o)}{\sum_{m=1}^{k} }[P_{i}^{J}(m)\eta_{i}(k)(1-\theta_{o}^{J}(m))]- \\
           &\sum_{m=0}^{k-1}\tilde{P}_{o}^{J}(m), \qquad \quad \text{if} \quad \varphi_{o}^{J}(k)\neq0 \\
           &0,\qquad \qquad \qquad \qquad  \enspace \text{if} \quad \varphi_{o}^{J}(k)=0\\
           &0, \qquad \qquad \qquad \qquad \enspace  \text{if} \quad  k=0
        \end{array}
    \right.
    \end{equation}

    \begin{equation}\label{eq:P2}
    \bar{P}_{i,o}^{J}(k) = \left\{
        \begin{array}{ll}
           \frac{\tilde{P}_{i,o}^{J}(k)}{\varphi_{o}^{J}(k)},  & \text{if} \quad \varphi_{o}^{J}(k)\neq0 \\
           0, & \text{if} \quad  \varphi_{o}^{J}(k)=0
        \end{array}
    \right.
    \end{equation}
From equations (\ref{eq:P1}) and (\ref{eq:P2}), we can know that the value of $\bar{P}_{i,o}^{J}(k)$ denotes the averaged number of pedestrians during current RED sub-sequence if and only if current traffic light $\theta_{o}(k)$ for stage $o$ is the last RED interval among current RED sub-sequence, otherwise, the value of $\bar{P}_{i,o}^{J}(k)$ is $0$.
\subsection{Solution algorithms}
\subsubsection{A MILP formulation for pedestrian total delay}
So far we have derived our pedestrian delay model as a linear cost function (\ref{eq:costped}) with its associated constraints which depict traffic light compatibility (\ref{eq:stageped}), volume dynamics (\ref{eq:dyn1a})-(\ref{eq:dyn1c}), capacity limitation (\ref{eq:capacityfun1}) or (\ref{eq:capacityfun2}) and hopping flow constraint (\ref{eq:hopflow}). Among all these constraints, (\ref{eq:stageped1b}), (\ref{eq:dyn1b}), (\ref{eq:dyn1c}), (\ref{eq:capacityfun1}), (\ref{eq:capacityfun2}) and (\ref{eq:hopflow}) can be converted into mixed integer linear constraints.\\

\indent Let $M$ chosen to be sufficiently big, e.g., $M \geq \text{max}_{k}(f_{ij}(k))$, then equation (\ref{eq:stageped1b}) can be translated into the following form:
       \begin{equation} \label{ineq1}
	       (\forall o \in \Omega, \forall (i,j) \in  h_{p}^{J}(o)) f_{ij}(k) \leq M*\theta_{o}(k) \\
	   \end{equation}
\begin{proposition}
Replacing Equation (\ref{eq:stageped1b}) with Inequalities (\ref{ineq1}) in the model leads to the same solution.
\end{proposition}

In order to decrease the number of decision variables, we substitute equations (\ref{eq:dyn1b}) and (\ref{eq:dyn1c}) into (\ref{eq:dyn1a}), then translate into the following form by applying the corresponding round-off measure:

       \begin{subequations}\label{eq:dyn3}
       \begin{align}
	       &P_{i}(k+1) - P_{i}(k) +  (\underset{i \in \mathcal{C}:(i,j)\in \cup_{J \in \mathcal{J}}\mathcal{F}_{p}^{J}}\sum{f_{ij}(k)})\triangle \leq \nonumber\\
           & I_{i}(k)\triangle+(1-\gamma_{i}(k))( \underset{i \in \mathcal{C}:(j,i)\in \cup_{J \in \mathcal{J}}\mathcal{F}_{p}^{J}}\sum{f_{ji}(k)})\triangle \label{eq:dyn3a}\\
           &P_{i}(k+1) - P_{i}(k) +  (\underset{i \in \mathcal{C}:(i,j)\in \cup_{J \in \mathcal{J}}\mathcal{F}_{p}^{J}}\sum{f_{ij}(k)})\triangle - \epsilon \geq  \nonumber\\
           &I_{i}(k)\triangle+(1-\gamma_{i}(k))( \underset{i \in \mathcal{C}:(j,i)\in \cup_{J \in \mathcal{J}}\mathcal{F}_{p}^{J}}\sum{f_{ji}(k)})\triangle -1 \label{eq:dyn3b}
       \end{align}
	   \end{subequations}

\begin{proposition}
Replacing Equation (\ref{eq:dyn1}) with Inequalities (\ref{eq:dyn3}) in the model leads to the same solution.
\end{proposition}

\indent Equation (\ref{eq:capacityfun1}) and (\ref{eq:capacityfun2}) can be translated into the following constraints by introducing binary variable $\delta_{o}(k)$:
       \begin{subequations} \label{condition}
       \begin{align}
	       & \delta_{o}(k) = 1 \Longleftrightarrow  k = 1
        \end{align}
	   \end{subequations}
\indent Then the logic relationship above can be rewritten as:
    	\begin{subequations}\label{eq:capselectcon}
        \begin{align}
        & M_{1}(\delta_{o}(k)-1) \leq k-1 \leq M_{1}(1 - \delta_{o}(k))\\
        & \hat{P}_{o}(k)-n(\triangle - I - \frac{L}{S_{p}}) \leq  M_{1}(1-\delta_{o}(k))\\
        & -\hat{P}_{o}(k)+n(\triangle - I - \frac{L}{S_{p}})-1\leq  M_{1}(1-\delta_{o}(k))-\epsilon\\
        & \hat{P}_{o}(k)-n(\triangle - I - \frac{L}{S_{p}}) \leq  M_{1}\theta_{o}(k-1)\\
        & -\hat{P}_{o}(k)+n(\triangle - I - \frac{L}{S_{p}})-1 \leq M_{1}\theta_{o}(k-1)-\epsilon\\
        & \hat{P}_{o}(k)-n\triangle \leq  M_{1}(1-\theta_{o}(k-1))\\
        & -\hat{P}_{o}(k)+n\triangle-1 \leq  M_{1}(1-\theta_{o}(k-1))-\epsilon
        \end{align}
        \end{subequations}
where $n = \frac{1}{0.27}$ if crosswalk width $W$ is smaller and equal to 3m, otherwise, $n = \frac{W}{0.81}$. $M_{1}$ is a sufficiently big integer, e.g., $M_{1} \geq \text{max}_{k}(H_p, \pm(\hat{P}_{o}(k)-n(\triangle - I - \frac{L}{S_{p}})), \pm(\hat{P}_{o}(k)-n\triangle))$, and $\epsilon$ is sufficiently small.
\begin{proposition}
Replacing Equation (\ref{eq:capacityfun1}) or (\ref{eq:capacityfun2}) with Inequalities (\ref{eq:capselectcon}) in the model leads to the same solution.
\end{proposition}

\indent Equation (\ref{eq:hopflow}) can be translated into the following form:
    	\begin{subequations}\label{eq:capselectcon}
        \begin{align}
        & (\forall o \in \Omega_{p}^{J}: (i,j) \in  h_{p}^{J}(o)) f_{ij}(k)\triangle \leq  \hat{P}_{o}(k)\\
        & (\forall o \in \Omega_{p}^{J}: (i,j) \in  h_{p}^{J}(o),\eta \in u_{p}^{J}(o)) f_{ij}(k)\triangle \leq  P_{i}(k)\eta_{i}(k)
        \end{align}
        \end{subequations}
\begin{proposition}
Replacing Equation (\ref{eq:hopflow}) with Inequalities (\ref{eq:capselectcon}) in the model leads to the same solution.
\end{proposition}

Finally, we can solve our pedestrian traffic scheduling problem as a mixed integer linear programming problem:

    \begin{equation}\label{eq:costped}
    \begin{aligned}
	\text{min} P_{D}&= \text{min}\sum_{J\in \mathcal{J}}\sum_{k=1}^{H_{p}}\sum_{i \in \mathcal{W}} [P_{i}^{J}(k)\\
    &- \underset{i \in \mathcal{C}:(i,j)\in \cup_{J \in \mathcal{J}}\mathcal{F}_{p}^{J}}\sum{f_{ij}(k)}]\triangle
    \end{aligned}
	\end{equation}
subject to:\\
        $ \underset{o \in \Omega_{p}^{J}}\sum\theta_{o}(k)=1, (\forall o \in \Omega_{p}^{J})(\forall k \in \mathbb{N}) \theta_{o}(k)\in \{0, 1\},\\$
        $(\forall o \in \Omega, \forall (i,j) \in  h_{p}^{J}(o)) f_{ij}(k) \leq M*\theta_{o}(k),\\$
        $P_{i}(k+1) - P_{i}(k) +  (\underset{i \in \mathcal{C}:(i,j)\in \cup_{J \in \mathcal{J}}\mathcal{F}_{p}^{J}}\sum{f_{ij}(k)})\triangle \leq \\$
        $ I_{i}(k)\triangle+(1-\gamma_{i}(k))( \underset{i \in \mathcal{C}:(j,i)\in \cup_{J \in \mathcal{J}}\mathcal{F}_{p}^{J}}\sum{f_{ji}(k)})\triangle,\\$
        $ P_{i}(k+1) - P_{i}(k) +  (\underset{i \in \mathcal{C}:(i,j)\in \cup_{J \in \mathcal{J}}\mathcal{F}_{p}^{J}}\sum{f_{ij}(k)})\triangle - \epsilon \geq  \\$
         $I_{i}(k)\triangle+(1-\gamma_{i}(k))( \underset{i \in \mathcal{C}:(j,i)\in \cup_{J \in \mathcal{J}}\mathcal{F}_{p}^{J}}\sum{f_{ji}(k)})\triangle -1,\\$
        $ M_{1}(\delta_{o}(k)-1) \leq k-1 \leq M_{1}(1 - \delta_{o}(k)),\\$
        $ \hat{P}_{o}(k)-n(\triangle - I - \frac{L}{S_{p}}) \leq  M_{1}(1-\delta_{o}(k)),\\$
        $ -\hat{P}_{o}(k)+n(\triangle - I - \frac{L}{S_{p}})-1\leq  M_{1}(1-\delta_{o}(k))-\epsilon,\\$
        $ \hat{P}_{o}(k)-n(\triangle - I - \frac{L}{S_{p}}) \leq  M_{1}\theta_{o}(k-1),\\$
        $ -\hat{P}_{o}(k)+n(\triangle - I - \frac{L}{S_{p}})-1 \leq M_{1}\theta_{o}(k-1)-\epsilon,\\$
        $ \hat{P}_{o}(k)-n\triangle \leq  M_{1}(1-\theta_{o}(k-1)),\\$
        $ -\hat{P}_{o}(k)+n\triangle-1 \leq  M_{1}(1-\theta_{o}(k-1))-\epsilon,\\$
        $(\forall o \in \Omega_{p}^{J}: (i,j) \in  h_{p}^{J}(o)) \quad f_{ij}(k)\triangle \leq  \hat{P}_{o}(k),\\$
        $(\forall o \in \Omega_{p}^{J}: (i,j) \in  h_{p}^{J}(o),\eta \in u_{p}^{J}(o)) \quad f_{ij}(k)\triangle \leq  P_{i}(k)\eta_{i}(k).$
\subsubsection{Discrete harmony search algorithm for pedestrian total delay and unhappiness} \label{sec:HS}
In the MILP formulation of the pedestrian delay problem stated above, an increase of the prediction horizon $H_{p}$ would exponentially increase the computational challenge. Additionally, the cost used in pedestrian unhappiness model is an nonlinear function, and the linear transformation would lead to the loss of subtle nonlinear characteristics of the model, meanwhile, directly convert piece-wise functions (\ref{eq:h})-(\ref{eq:P2}) into mixed integer linear constraints will make the problem even more computationally expensive. In order to overcome the above deficiencies of mathematical techniques, many heuristic optimization approaches have been developed. Although the final result may not be the global optima, a good signal setting (near optimum) in real time is more urgently needed rather than the best solution with costly time. \\
\indent The harmony search (HS) algorithm, mimicking the musical performance process to reach a perfect state of harmony, is one of the meta-heuristic algorithms \cite{geem2001new}. In the HS, four parameters play the major roles in obtaining the best 'harmony': Harmony Memory Size (HMS) , Harmony Memory Considering Rate (HMCR), Pitch Adjustment Rate (PAR) and the number of improvisations (NI). HMS defines the number of existed solution vectors in the Harmony Memory (HM), HMCR is the rate of selecting the values from the HM, PAR sets the adjustment rate of the pitch from the HM and NI is the number of iterations. In discrete harmony search (DHS) algorithm, the solution $X_{i}$ is the $i^{th}$ solution in the HM, and it is a $n$ dimensional vector $X_{i} = \{X_{i}(1), X_{i}(2), \cdots, X_{i}(n)\}$, which consists of $n$ decision variables of the optimization problem. The specific steps of the DHS algorithm can be given as:

Step 1 - Initialization

The element $X_{i}(k)$ of each harmony vector $X_{i}$ in the HM can be generated as follows:
  \begin{equation}\label{eq:HS1}
  \begin{aligned}
    &X_{i}(k) = LB(k) + (UB(k)-LB(k))*rand(),\\
    & \forall k=1,2,..,n \quad \text{and} \quad \forall i = 1,2,...HMS
  \end{aligned}
  \end{equation}
  where $LB(k)$ and $UB(k)$ are the lower bound and upper bound of the decision variable $X_{i}(k)$, respectively, and $rand()$ randomly generates value between $0$ and $1$.

Step 2 - Improvisation

The new harmony vector $X_{new}$ is created by setting the HS parameters: HMCR, PAR and random selection. Firstly, if a random value $rand()$ in the range of $[0, 1]$ is less than HMCR, the decision variable $X_{new}(k)$ is randomly selected from any value in the HM, otherwise, $X_{new}(k)$ is generated according to the upper and lower bounds of the decision variable $X_{i}(k)$. After improvisation, the pitch adjustment with a probability of PAR is applied to the new harmony vector if this vector is selected from the current HM. The whole process is illustrated as follows:
\begin{algorithm}
\caption{Generating new harmony vector}
	If ($rand1()<HMCR$)
	\begin{eqnarray}
	X_{new} (k)=X_{i} (k)
	\end{eqnarray}
	\begin{itemize}
		\item[] If ($rand2()<PAR$)
			\begin{eqnarray}
			X_{new} (k)=X_{new} (k)\pm rand()\times BW
			\end{eqnarray}
	\end{itemize}
	Else
	\begin{eqnarray}
	X_{new} (k)=LB(j)+(UB(k)-LB(k))\times rand()
	\end{eqnarray}
\end{algorithm}

 where $BW$ is the bandwidth, which is defined to find possible neighborhood solutions, and $i\in \{1,2, \cdots,\text{HMS}\}$, $X_i (j)$ is randomly selected from HM.

Step 3 - Update and Stopping condition

After determining the objective function values in HM, the harmony vector which give the worst objective function value is removed from the HM. Then checking the termination criterion, e.g, the number of NI, if the predetermined criterion is met, the searching process is terminated, otherwise, go back to Step 2.

\subsection{Numerical examples and sensitivity analysis} \label{sec:2-experiment}
In this section we firstly consider a pedestrian delay model, then a pedestrian unhappiness model is discussed at the second part. In order to conveniently coordinate with the vehicle network in later section, both models here are considered at a network-based level.
The network has $n_{h}$ horizontal links per row and $n_{v}$ vertical links per column, and the total number
of junctions is $n_{J} = n_{h}n_{v}$. Each intersection only has two stages: pedestrians cross the road horizontally and vertically. The information
required by solving the pedestrian traffic scheduling problem at time $k$ includes the current pedestrian demand $P_{i}(k)$, pedestrian diversion ratio $\alpha_{i}(k)$ and $\beta_{i}(k)$, pedestrian departure ratio $\gamma_{i}(k)$, and a prediction of incoming flow $I_{i}(k)$ to $I_{i}(k+N)$ ($N$ is the prediction step) for $i \in \{1,2,3,4\}$ and $J \in \mathcal{J}$, which can be measured by sensors. Considering the higher computational complexity from unhappiness model, the MILP formulation of traffic light scheduling is only utilized for pedestrian delay model, then the discrete harmony search algorithm is used for both delay and unhappiness models. The parameters used in case studies are summarized below:

    \begin{table}[!ht]
		\renewcommand{\arraystretch}{0.9}
		\caption{Parameters used in case studies}
       \label{label11}
		\centering
          \begin{tabular}{l|l|l}
			\hline
			\multicolumn{1}{c|}{\textbf{Parameters}}
          & \multicolumn{1}{c|}{\textbf{Descriptions}} & \multicolumn{1}{c}{\textbf{Associated Values}} \\ \hline
           $\triangle$ & sampling period & 15s\\
           $S_{p}$ & pedestrian average speed & 1.2m/s\\
           $L$ & crosswalk length & 8.5m\\
           $W$ & crosswalk width & 4m\\
           $I$ & start-up time & 3.2s\\
\hline
		\end{tabular}
	\end{table}

\subsubsection{pedestrian delay model}
The urban pedestrian traffic signal scheduling problem with consideration of minimizing delay times can be solved by an optimization solver to obtain an optimal traffic control signal strategy. During the experiment, a model predictive control strategy is adopted, namely, after an optimal traffic signal profile is calculated over a finite time horizon based on the current measurement of information from sensors, only the first interval of optimal profile is implemented, and this process continues till the end of the experiment. GUROBI \cite{optimization2012gurobi} is used in MATLAB to solve the problems on a PC with an Intel(R) Core(TM) i7-3770 CPU @3.40GHz and RAM 8GB.\\
\indent The Table \ref{table_ped1} illustrates the potential complexity involved in solving this problem, which lists the total number of decision variables $N_{x}$ and the total number of constraints $N_{cons}$ with $N = 5$ and $N = 10$ under different network sizes. Clearly, we need to solve a large mixed integer linear programming problem for a large urban traffic network, so some efficient computational algorithms need to be considered once large network involved in order to ensure proper real-time decisions.

\begin{table}[!ht]
\centering 
\caption{NUMBERS OF DECISION VARIABLES AND CONSTRAINTS}
\renewcommand{\arraystretch}{1.0}
\label{table_ped1}
\begin{tabular}{|l|l|l|l|l|}
\hline
\hline
$n_{v}=n_{h}$& \multicolumn{2}{c|}{N=5} & \multicolumn{2}{c|}{N=10} \\ \hline
      & $N_{x}$     & $N_{cons}$ & $N_{x}$     & $N_{cons}$  \\ \hline
3     & 1710        & 3843       & 3420        & 7686        \\ \hline
4     & 3040        & 6832       & 6080        & 13664       \\ \hline
5     & 4750        & 10675      & 9500        & 21350       \\ \hline
6     & 6840        & 15372      & 13680       & 30744       \\ \hline
7     & 9310        & 20923      & 18620       & 41864       \\ \hline
8     & 12160       & 27328      & 24320       & 54656       \\ \hline
9     & 15390       & 34587      & 30780       & 69174       \\ \hline
10    & 19000       & 42700      & 38000       & 85400       \\ \hline
\end{tabular}
\end{table}

We apply our mixed integer linear programming to the network with different sizes and different time horizons $H_{p}=N\triangle$ to check how fast the scheduling problem can be solved, and we tested cases of $n_{v}=n_{h}=3, 4, 5, 6, 7, 8, 9, 10$ and $H_{p}=15, 30, 45, 60, 75$ seconds, respectively. Table \ref{table_ped2} depicts the processing time of the pedestrian model associated with Table \ref{table_ped1}. In Table \ref{table_ped2}, we can see that GUROBI can solve the 10$\times$10 junctions with a 30-second scheduling horizon in only $1.4040$ second, which is equivalent to solving a mixed integer linear programming problem with 8200 integer and binary decision variables and 17080 constraints. This is definitely workable in practical real-time scheduling environment, because only $9.3600\%$ ($1.4040/15$) of the time is used for computation. However, when the network and prediction horizon become larger, the probability of successfully solving the scheduling problem in real time becomes lower. For example, a 10$\times$10 network with a 90-second prediction horizon needs $59.9582$ second to obtain a (near-)optimal signal profile, which has already reached nearly 4 times of the sampling period $15$ second, so computational challenge requires more efficient algorithms.

\begin{table}[!ht]
\caption{PROCESSING TIME OF PEDESTRIAN DELAY MODEL + MILP}
\renewcommand{\arraystretch}{1.0}
\label{table_ped2}
\begin{center}
\resizebox{\columnwidth}{!}{
\begin{tabular}{|c|c|c|c|c|c|}
\hline
\hline
 $n_{v} = n_{h}$ & $H_{p}=$30s & $H_{p}=$45s &$H_{p}=$60s & $H_{p}=$75s & $H_{p}=$90s\\
\hline
$3$ & 0.2971 &0.4853 &1.0638 &2.5521&4.8843\\
\hline
$4$ & 0.3377 &0.8039 &1.3923 &3.8753&9.0773\\
\hline
$5$ &0.4152 &1.0774&2.4544 &6.4511&12.2866\\
\hline
$6$ & 0.5186 &1.4292 &3.3435 &9.1050&19.6239\\
\hline
$7$ &0.6221 &1.9756&4.5597 &12.2719&27.3518\\
\hline
$8$ &0.8919 &2.8823&6.7923 &19.0293&44.5946\\
\hline
$9$ &1.1035 &3.7512&8.5685 &22.1222&50.8451\\
\hline
$10$ &1.4040 &4.7606&11.9072 &28.1614&59.9582\\
\hline
\end{tabular}
}
\end{center}
\end{table}

Considering the computational challenge when network and prediction horizon become larger, a discrete harmony search algorithm discussed in Section \ref{sec:HS} is used to solve the problem. The population size (Harmony Memory Size) is 1000, and the maximum number of generations is 1000. Parameters HMCR and PAR are set as 0.95 and 0.5, respectively. The sizes of the traffic grids are from 3$\times$3 to 10$\times$10, and the prediction horizons are 30 seconds, 60 seconds and 90 seconds, respectively.

The comparison results between MILP and DHS are shown in Table \ref{MILP:DHS}. By applying the optimal traffic signal settings obtained from MILP to the entire network, the actual delay times $J$ are listed in the table, and the computed delay times $J^{*}$ by applying near-optimal signals from DHS are also listed in the table. Additionally, we can see that the difference degree $\frac{J-J^{*}}{J^{*}}$  between delay times obtained from two different algorithms is very close. Although the gap between $J$ and $J^{*}$ become larger when both network size and prediction horizon increase, the difference degree is still below 20$\%$ for a 10$\times$10 network with a 90-second prediction horizon. Certainly, the outcome quality will gradually decrease with the growth of the network size and prediction horizon. However, the computational times, refer to the second column of Table \ref{MILP:DHS}, significantly drop when applying DHS algorithm, especially for the large traffic network with longer prediction horizon. For example, a 10$\times$10 traffic system with 90s prediction horizon require 59.9582s in Table \ref{table_ped2} to obtain the optimal solution by MILP. However, with the same scale setting, the DHS algorithm only needs 0.6020s to get a near-optimal solution, which has approximately 20$\%$ gap of the optimal solution and is absolutely feasible in realistic traffic assignment. Furthermore, the gap can even be neglected when a 30s prediction horizon is conducted.

\begin{table}[!ht]
\centering
\caption{A COMPARISON BETWEEN COMPUTED AND ACTUAL TIME DELAY FOR
DHS WITH MILP}
\label{MILP:DHS}
\begin{tabular}{|c|c|c|c|c|}
\hline
\hline
Instance & \multicolumn{4}{c|}{$H_{p}=$30s}                  \\ \hline
$n_{v} = n_{h}$ & Processing Time & $J$     & $J^{*}$ & $\frac{J-J^{*}}{J^{*}}$ \\ \hline
3        & 0.0230          & 2713   & 2713  & 0.00\%    \\ \hline
4        & 0.0370          & 4572   & 4566  & 0.13\%    \\ \hline
5        & 0.0540          & 6697   & 6654  & 0.65\%    \\ \hline
6        & 0.0780          & 10006  & 9949  & 0.57\%    \\ \hline
7        & 0.1030          & 13296  & 13023 & 2.09\%    \\ \hline
8        & 0.1320          & 17279  & 16880 & 2.36\%    \\ \hline
9        & 0.1640          & 21533  & 20947 & 2.79\%    \\ \hline
10       & 0.2010          & 27272  & 26157 & 4.26\%    \\ \hline
\hline
Instance & \multicolumn{4}{c|}{$H_{p}=$60s}                  \\ \hline
$n_{v} = n_{h}$& Processing Time & $J$      & $J^{*}$    & $\frac{J-J^{*}}{J^{*}}$ \\ \hline
3        & 0.04            & 5822   & 5730  & 1.61\%    \\ \hline
4        & 0.068           & 9731   & 9653  & 0.80\%    \\ \hline
5        & 0.103           & 14828  & 14469 & 2.48\%    \\ \hline
6        & 0.151           & 22161  & 21025 & 5.40\%    \\ \hline
7        & 0.203           & 30001  & 27999 & 7.15\%    \\ \hline
8        & 0.26            & 39025  & 36326 & 7.43\%    \\ \hline
9        & 0.327           & 49375  & 44710 & 10.43\%   \\ \hline
10       & 0.401           & 62038  & 55965 & 10.85\%   \\ \hline
\hline
Instance & \multicolumn{4}{c|}{$H_{p}=$90s}                  \\ \hline
$n_{v} = n_{h}$& Processing Time & $J$      & $J^{*}$   & $\frac{J-J^{*}}{J^{*}}$ \\ \hline
3        & 0.0590          & 9574   & 9348  & 2.42\%    \\ \hline
4        & 0.1000          & 16134  & 15645 & 3.13\%    \\ \hline
5        & 0.1550          & 24208  & 22497 & 7.61\%    \\ \hline
6        & 0.2280          & 36853  & 33507 & 9.99\%    \\ \hline
7        & 0.3050          & 50888  & 44130 & 15.31\%   \\ \hline
8        & 0.3910          & 67134  & 59545 & 12.74\%   \\ \hline
9        & 0.4960          & 85724  & 72350 & 18.49\%   \\ \hline
10       & 0.6020          & 106826 & 90351 & 18.23\%   \\ \hline
\end{tabular}
\end{table}

\subsubsection{pedestrian unhappiness model}

Considering the computational complexity from pedestrian unhappiness model, the DHS algorithm and the new formulation of pedestrian unhappiness discussed in Section \ref{sec:unhappy2} are directly adopted. And the parameters used in DHS algorithm are the same as assignments in pedestrian delay model, namely, HMS is 1000, and the maximum number of generations is 1000. Parameters HMCR and PAR are set as 0.95 and 0.5, respectively.

In order to show that unhappiness model do provide fairness to those few pedestrians with the same forward direction, the switching frequency of traffic signals between two models are listed in Table \ref{Ped_Freq}. The switching frequency (SF) of traffic signals is defined as the portion of the changed signals in next time interval of entire network with respect to the total number of traffic signals of entire network, and if the prediction step $N$ is more than 2, then corresponding SF shown in Table \ref{Ped_Freq} is the averaged value of the switching frequencies with corresponding quantity $N-1$. From Table \ref{Ped_Freq}, it is obvious to see that the values of SF are all very small in delay model, on the contrary, the values of SF are quite large in unhappiness model. The reason why traffic signals switch infrequently in pedestrian delay model is that the delay times are linearly increase, and the capacity for each GREEN time interval is not a fixed value, which depends on the position of current GREEN time interval among the current GREEN sub-sequence (the capacity for the first GREEN time interval in GREEN sub-sequence is a small value, and the capacity for other GREEN time intervals in GREEN sub-sequence is a large value), so a longer GREEN time will be applied in one direction in order to maximize the performance. While in pedestrian unhappiness model, the exponential increase of the pedestrian unhappiness will compel the switch of the signals, and giving any direction more GREEN time will pay even greater cost on the other direction, namely, the privilege obtained by keeping a longer GREEN time in one direction to obtain more capacity is no longer comparable with the the exponential increase of the pedestrian unhappiness in the other direction. Accordingly, the phenomenon frequently occur in pedestrian delay model that few pedestrians wait for a long time will be significantly reduced under this framework.

\begin{table}[!ht]
\centering
\caption{COMPARISON OF SWITCHING FREQUENCY OF TRAFFIC SIGNALS BETWEEN PEDESTRIAN DELAY MODEL AND PEDESTRIAN UNHAPPINESS MODEL}
\label{Ped_Freq}
\resizebox{\columnwidth}{!}{
\begin{tabular}{|c|c|c|}
\hline
\hline
Instance & \multicolumn{2}{c|}{$H_{p}=$30s}                           \\ \hline
$n_{v} = n_{h}$& Pedestrian Delay Model & Pedestrian Unhappiness Model \\ \hline
3        & 0.0000                 & 0.8888                       \\ \hline
4        & 0.0000                 & 0.9375                       \\ \hline
5        & 0.0000                 & 0.9200                       \\ \hline
6        & 0.0000                 & 0.8611                       \\ \hline
7        & 0.0204                 & 0.8775                      \\ \hline
8        & 0.0312                 & 0.8594                       \\ \hline
9        & 0.0123                 & 0.6914                      \\ \hline
10       & 0.0300                 & 0.8000                       \\ \hline
\hline
Instance & \multicolumn{2}{c|}{$H_{p}=$60s}                           \\ \hline
$n_{v} = n_{h}$& Pedestrian Delay Model & Pedestrian Unhappiness Model \\ \hline
3        & 0.1111                 & 0.7407                       \\ \hline
4        & 0.0833                 & 0.6875                       \\ \hline
5        & 0.1067                 & 0.6800                       \\ \hline
6        & 0.1111                 & 0.7500                       \\ \hline
7        & 0.1428                 & 0.7755                       \\ \hline
8        & 0.1354                 & 0.7240                       \\ \hline
9        & 0.1564                 & 0.7078                       \\ \hline
10       & 0.1900                 & 0.7200                       \\ \hline
\hline
Instance & \multicolumn{2}{c|}{$H_{p}=$90s}                           \\ \hline
$n_{v} = n_{h}$& Pedestrian Delay Model & Pedestrian Unhappiness Model \\ \hline
3        & 0.0667                 & 0.5555                       \\ \hline
4        & 0.1125                 & 0.5375                       \\ \hline
5        & 0.1440                 & 0.6528                       \\ \hline
6        & 0.1556                 & 0.6202                       \\ \hline
7        & 0.1716                 & 0.5037                       \\ \hline
8        & 0.2031                 & 0.6625                       \\ \hline
9        & 0.2247                 & 0.6351                      \\ \hline
10       & 0.2340                 & 0.6608                       \\ \hline
\end{tabular}
 }
\end{table}

\section{Formulation of A Traffic Light Scheduling Problem for pedestrian-vehicle mixed-flow network}\label{sec:3}
\subsection{A vehicle flow model}

Unlike the uncontinuity of a pedestrian traffic network, the vehicle traffic network is connected by a set of road links and junctions. We introduce a discrete-time vehicle model based on an extended cell transmission flow dynamic model. Similar to the pedestrian model, $\mathcal{L}$ denotes the set of links ($L\in \mathcal{L}$) and $ \mathcal{J}$ denotes the set of junctions ($J \in \mathcal{J}$), let $\Omega_{v}^{J}$ be the set of vehicle stages in junction $J$, $\mathcal{F}_{v}^{J}$ the set of all vehicle streams in junction $J$: $\mathcal{F}_{v}^{J}$ $\in$ $\mathcal{L} \times \mathcal{L}$, and $h_{v}^{J}$ the association of each stage to relevant compatible streams, $h_{v}^{J}: \Omega_{v}^{J} \to 2^{\mathcal{F}_{v}^{J}}$.
The following assumptions are considered in order to fit this deterministic model.
	\begin{enumerate}[label=\textbf{\arabic*},start=1]
			\item The network boundary demand for entrance and exit are known.
            \item No traffic demand is generated inside the network.	
        	\item The link turning ratios of the network are known.
            \item Each vehicle delayed only by traffic lights will finally leave the network.
		\end{enumerate}
Assumption 2 can be easily relaxed by slightly changing the volume dynamic constraints in equation (\ref{eq:ldyn1}) with the consideration of the demand and the exit flow within the link, and we give this assumption here to better serve our case studies in Section \ref{sec:3}.

\subsubsection{Parameters}
	 To better describe the vehicle model, the parameters used in later part are summarized in TABLE \ref{label12}.
    \begin{table}[!ht]
		\caption{Notations for the vehicle model}
       \label{label12}
		\centering
          \begin{tabular}{p{0.22\linewidth}p{0.7\linewidth}}
		  \hline
		  \textbf{Parameter} & \textbf{Definition}\\
		  \hline
            $w$ &  A set of signal stages to control antagonistic vehicle flow streams, $w \in \Omega_{v}^{J}$.\\
			$\theta_{w}(k)$  &  The vehicle traffic light in associated stage $w$. \\
            $f_{ij}(k)$ &  The flow rate from link $i$ to link $j$ during interval $k$. \\
			$C_{i}(k)$ & The number of  vehicles in link $i$ during interval $k$.\\
			$\hat{C}_{i}$ &  The maximum volume of link $i$.\\

			$d_{j}(k)$ & The number of incoming vehicles of link $j$ at time interval $k$.\\
			
			$s_{j}(k)$ & The number of outgoing vehicles of link $j$ at time interval $k$.\\
			
            $\lambda_{ij}(k)$ & Turning ratio of vehicles in link $i$ towards link $j$ during time interval $k$,
            $\sum\limits_{j\in\mathcal{L}:(i,j)\in\cup_{J\in\mathcal{J}}\mathcal{F}_{v}^{J}} \lambda_{ij}(k)=1$.\\

\hline
		\end{tabular}
	\end{table}
\subsubsection{Stage Constraints}
All associated flows with stage $w$ are zero if the stage traffic light is RED and only one active stage exists for one junction $J$ in each time interval $k$, and $\mathbb{N}$ denotes the set of natural numbers, which can be described as:
    	\begin{subequations}\label{eq:conn1}
        \begin{align}
	       &(\forall w \in \Omega_{v}^{J}) \theta_{w}(k)=0 \nonumber \\
           &\Rightarrow (\forall (i,j) \in  h_{v}^{J}(w)) f_{ij}(k)=0 \label{eq:stage1a}\\
           &\underset{w \in \Omega_{v}^{J}}\sum\theta_{w}(k)=1 \label{eq:stage1b}\\
           &(\forall w \in \Omega_{v}^{J})(\forall k \in \mathbb{N}) \theta_{w}(k)\in \{0, 1\} \label{eq:stage1c}
        \end{align}
	    \end{subequations}

\subsubsection{Link Dynamics and Flow Dynamics}
According to the conservation of vehicles, the dynamic of each link follows the rule:
    	\begin{equation}\label{eq:ldyn1}
	       (\forall k \in \mathbb{N})C_{j}(k+1) = C_{j}(k) + \triangle(d_{j}(k)-s_{j}(k))
	    \end{equation}
    \begin{equation}\label{eq:fdyn1}
    \begin{aligned}
	f_{ij}(k)= & \text{min} \{ \lfloor {\lambda_{ij}(k)C_{i}(k)}\rfloor, \hat{C}_{j}(k)- C_{j}(k)+ s_{j}(k), \\
               & \lfloor{l_{ij}(k)v_{i}^{*}d^{*}\triangle} \rfloor \}
    \end{aligned}
	\end{equation}	

Incoming flow $d_{j}(k)$ and outgoing flow $s_{j}(k)$ can be captured by the following equations:

    	\begin{subequations}\label{eq:descrp1}
        \begin{align}
	       &d_{j}(k) = \underset{i \in \mathcal{L}:(i,j)\in \cup_{J \in \mathcal{J}}\mathcal{F}_{v}^{J}}\sum{f_{ij}(k)} \label{eq:descrp1a}\\
           &s_{j}(k) = \underset{i \in \mathcal{L}:(j,i)\in \cup_{J \in \mathcal{J}}\mathcal{F}_{v}^{J}}\sum{f_{ji}(k)}\label{eq:descrp1b}
        \end{align}
	    \end{subequations}

The flow dynamics constraint described in Condition (\ref{eq:fdyn1}) indicate that the outgoing flow $f_{ij}(k)$ is determined by the current upstream volume $\lambda_{ij}(k)C_{i}(k)$ of link $i$, the downstream remaining space $\hat{C}_{j}(k)- C_{j}(k)+ s_{j}(k)$, and the capacity from link $i$ to link $j$ in one time interval, which is captured by the critical speed $v^{*}$, critical density $d^{*}$, time interval $\triangle$ and the speed category $l_{ij}(k)$. We construct the outflow $f_{ij}(k)$ within each time interval as a nonlinear mixed logical switching function over the upstream link's density, the downstream link's density and capacity, and the drivers' psychological behavior respond to the past continuous state of traffic lights.

The intention of the psychological response is that drivers are more likely to keep a high speed $v_{ij}(k)$ if the stage $w$ has been active for the past $r$ intervals and the downstream link has sufficient space to receive the flow. Thus the corresponding speed is $v_{ij}(k) = l_{ij}(k)v^{*}$, where the speed category $l_{ij}(k)$ can be approximated by a set of discrete speed levels: $l_{ij}^{0} \geq ... \geq l_{ij}^{r} > 0$ as follows:

    	\begin{subequations}\label{eq:speedca1}
        \begin{align}
	       & l_{ij}(k) = \sum\limits_{p = 0}^{r} \delta_{ij}^{p}(k)l_{ij}^{p} \label{eq:speedca1a} \\
           & \sum\limits_{p = 0}^{r} \delta_{ij}^{p}(k)l_{ij}^{p} - \theta_{w}(k) = 0 \label{eq:speedca1b}\\
           &(\forall q \in [0,r-1])(1-\theta_{w}(k-q-1))\prod\limits_{p=0}^{q}\theta_{w}(k-p)=1 \nonumber \\
           & \qquad \qquad \Leftrightarrow \delta_{ij}^{r-q}(k)=1 \label{eq:speedca1c}\\
           & \prod \limits_{p = 0}^{r} \theta_{w}(k-p)=1 \Leftrightarrow \delta_{ij}^{0}(k)=1 \label{eq:speedca1d}\\
           &(\forall p \in [0,r]) \delta_{ij}^{p}(k) \in \{0,1\} \label{eq:speedca1e}
        \end{align}
	    \end{subequations}

Condition (\ref{eq:speedca1b}) indicates that all $\delta_{ij}^{p}(k) = 0$ if the traffic light of stage $w$ is RED, then from Condition (\ref{eq:speedca1a}), $l_{ij}(k)=0$. If the traffic light of stage $w$ is GREEN, then $l_{ij}(k)$ can only choose one speed level $l_{ij}^p$ determined by the following Conditions (\ref{eq:speedca1c})-(\ref{eq:speedca1d}), which indicate that the actual speed can be obtained when looking back the number of consecutive green light intervals from current interval $k$. The larger the number of consecutive green light intervals, the higher the speed category.
\subsubsection{Optimal Control for Vehicle Delay}
Our cost function to minimize the total delay time for the whole network within $N$ time intervals can be formulated as follows:
    \begin{equation}\label{eq:costveh}
    \begin{aligned}
	\text{min}& V_{D}=\\
      &\text{min} \sum_{i \in \mathcal{L}}\sum_{k=1}^{N} [C_{i}(k)-
      \frac{L_{i}}{v_{i,max}}\sum \limits_{j\in \mathcal{L}:(i,j)\in \cup_{J\in\mathcal{J}}\mathcal{F}_{v}^{J}}f_{ij}(k)]\triangle
    \end{aligned}
	\end{equation}
Where $v_{i, max}$ is the free flow speed of link $i$, $L_{i}$ is the length of link $i$, and $\triangle$ is the sampling interval period.
\subsection{Combination of pedestrian model and vehicle model} \label{subsecintegration}
The two delay models that we introduced above are both discrete-time models and they focus on different objects - pedestrians and vehicles. Our goal is to see how the vehicle traffic total delay and its corresponding traffic lights will be affected if we assign different priorities to the pedestrian side. For this purpose, we propose the following formulation based on weighted sum. To ensure that the weights assigned to the pedestrians and vehicles are in the same magnitude, we need the following scale-down operation.
\subsubsection{Function Transformations}
The approach we use to scale down objective functions is given as follow:
    \begin{equation}\label{eq:transPed}
	P_{D}^{trans} = \frac{P_{D}}{|P_{D}^{max}|}
	\end{equation}
    \begin{equation}\label{eq:transVeh}
	V_{D}^{trans} = \frac{V_{D}}{|V_{D}^{max}|}
	\end{equation}
where $P_{D}$ and $V_{D}$ are our original pedestrian and vehicle cost function respectively, while $P_{D}^{max}$ and $V_{D}^{max}$ are the absolute maximum of $P_{D}$ and $V_{D}$ respectively, which can be obtained by changing the direction of the objective functions ($\ref{eq:costped}$) and ($\ref{eq:costveh}$) from $\text{min}$ to $\text{max}$. This approach can effectively make the magnitudes of two different-scale problems comparable.
\subsubsection{Weighted Sum Method}
One of the common scalarization methods for multi-objective optimization is the weighted sum method in which all objective functions are combined to form a single cost function. Since we want to see how pedestrians influence the vehicle, the corresponding weighted function is chosen as follow:
    \begin{equation} \label{eq:integrate}
	U_{D} = V_{D}^{trans} + mP_{D}^{trans}
	\end{equation}
where $m$ is associated weighting coefficients.
\subsubsection{Constraints between vehicles and pedestrians}
Considering the conflict between vehicles and pedestrians, we give the constraints below to avoid collision:

%

       \begin{equation}\label{eq:combine1}
       \begin{aligned}
	      \underset{n \in \Omega_{p+v}^{J}} \sum{\theta_{n}(k)} = (or \leq)1
        \end{aligned}
	   \end{equation}
Where $\Omega_{p+v}^{J}$ denotes the total stages to control all antagonistic flow streams for both vehicles and pedestrians in junction $J$, and $n \in \Omega_{p+v}^{J}$.
The objective of our final integrated model is captured by equation (\ref{eq:integrate}) and its associated constraints are represented as: (\ref{eq:stageped})-(\ref{eq:hopflow}), (\ref{eq:conn1})-(\ref{eq:speedca1}) and (\ref{eq:combine1}).

\subsection{Case study - A one-way vehicle network with combination of a pedestrian network}
\begin{figure}[!ht]
			\centering
			\includegraphics[width=2.2in]{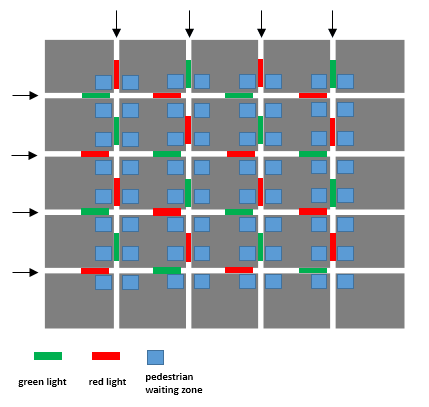}
			\caption{A simplified integrated urban traffic network}
			\label{(fig_sim4)}
		\end{figure}

Similar to the pedestrian network introduced in Section \ref{sec:2-experiment}, the integrated network also have $n_{h}$ horizontal links per row and $n_{v}$ vertical links per column. We now introduce a simplified vehicle traffic network which only has traffic flows moving from left to right horizontally and from
top to bottom vertically, accordingly, each junction only have two stages. An incoming boundary flow means a flow comes from outside the network to a boundary link which
locates either at the first row or at the first column of the network shown in Fig. \ref{(fig_sim4)}.
The information required by solving the vehicle traffic scheduling problem at time $k$ includes the incoming boundary flows $f_{ij}(k)$ to $f_{ij}(k+N)$, the current link volume $C_{i}(k)$, the past traffic light state $\theta_{w}(k)$. Also, we only consider
two speed categories and associated speed levels are $l_{ij}^0=0.3$ and $l_{ij}^0=0.15$, the maximal volume $\hat{C_{j}}$ for each link $j$ is 100 and $L_{i}/v_{i,max}=1$ is
chosen in simulation. And the information required by the pedestrian network is the same as the individual pedestrian network discussed in Section \ref{sec:2-experiment}.\\
\indent Then we combine the pedestrian network with this simplified vehicle network to see how pedestrians influence vehicles, which is shown in Fig. \ref{(fig_sim4)}. In order to reflect the more real effects from the pedestrian side to the entire integrated network, GUROBI \cite{optimization2012gurobi} is also used in MATLAB to obtain the optimal solution.
According to the transformation method and weighted sum method introduced in Section \ref{sec:2}, we generate the results which illustrate how the vehicle traffic light changes when giving different weights on the pedestrian side. Table \ref{integrate1} and Table \ref{integrate2} show the state of the vehicle traffic light in 3$\times$3 and 4$\times$4 networks with a 30s prediction horizon respectively, $0$ means the traffic light is RED and the corresponding horizontal flow is blocked, and $1$ means the traffic light is GREEN and the corresponding vertical flow is blocked. The weights that we give in two Tables below are all the turning points, which means these weights are exactly associated with the change of one or several vehicle traffic lights, also, the specific change state of each traffic light is highlighted in Tables below. The prediction information, e.g., the incoming boundary flows $f_{ij}(k)$ to $f_{ij}(k+N)$, will become increasingly inaccurate with the increase of the prediction horizon, which will also lead to the lower accuracy of the influence bring from pedestrians to vehicles, so we choose 30s prediction horizon in our simulation. The switching frequency (SF) is defined as the number of the changed traffic light at the turning weight with respect to the traffic lights of the entire network and the entire prediction horizon, which is different with the switching frequency discussed in Section \ref{sec:2}. SF in Section \ref{sec:2} compare the traffic light states at two connected prediction horizon, but SF here compare the traffic light states at associated turning weight. And the specific values of SF are listed in Table \ref{sf1} and Table \ref{sf2}, respectively, which are corresponding to the Table \ref{integrate1} and Table \ref{integrate2}.

\begin{figure}[!ht]
		\centering
		\includegraphics[width=0.9\linewidth]{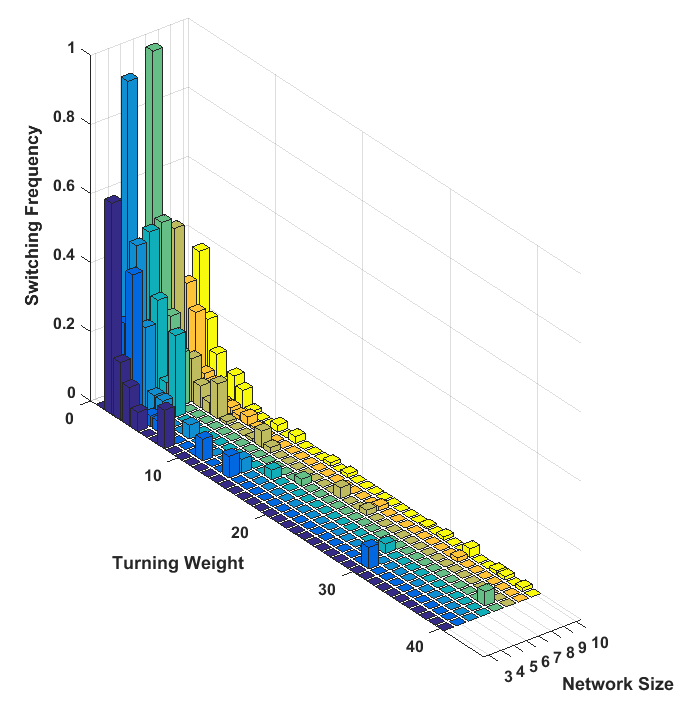}
	\caption{Switching frequency of traffic signals at turning weight from network 3$\times$3 to network 10$\times$10}
    \label{fig:sf3}
\end{figure}

\begin{figure}[!ht]
	\centering
		\includegraphics[width=0.9\linewidth]{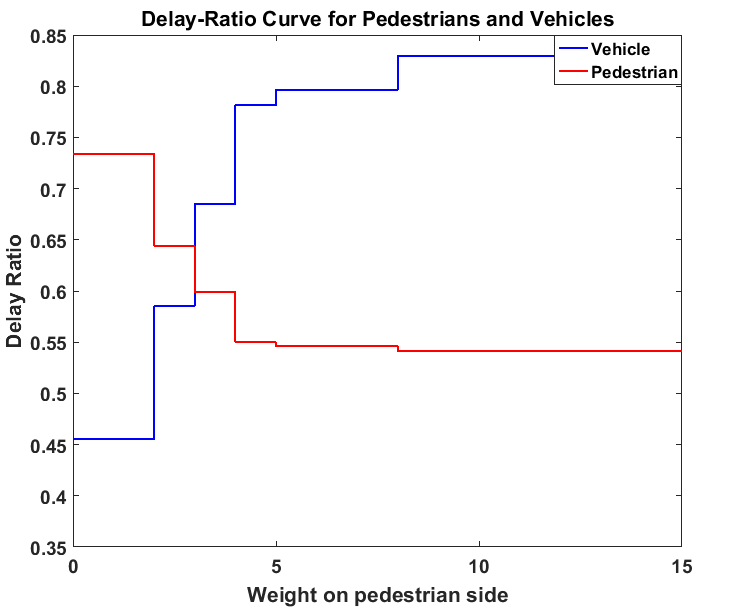}
		\caption{The change of delay-ratio curve under different weight of traffic composition from network 3$\times$3}
    \label{fig:compare1}
\end{figure}

\begin{figure}[!ht]
	\centering
	\includegraphics[width=0.9\linewidth]{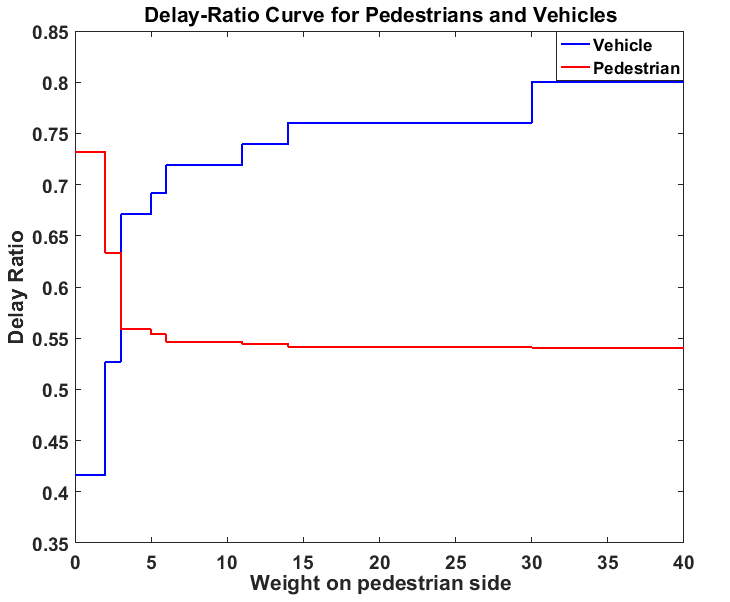}
	\caption{The change of delay-ratio curve under different weight of traffic composition from network 4$\times$4}
\label{fig:compare2}
\end{figure}

\begin{figure}[!ht]
	\centering
	\includegraphics[width=0.9\linewidth]{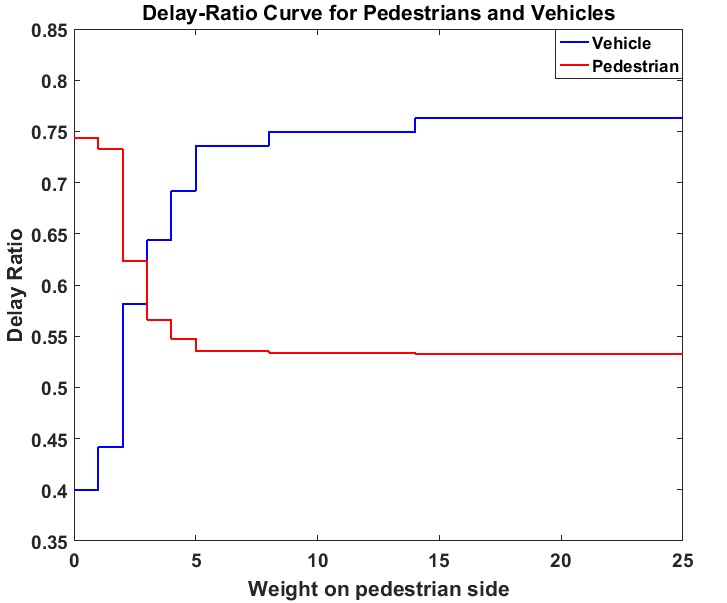}
	\caption{The change of delay-ratio curve under different weight of traffic composition from network 5$\times$5}
\label{fig:compare3}
\end{figure}

\begin{figure}[!ht]
	\centering
	\includegraphics[width=0.9\linewidth]{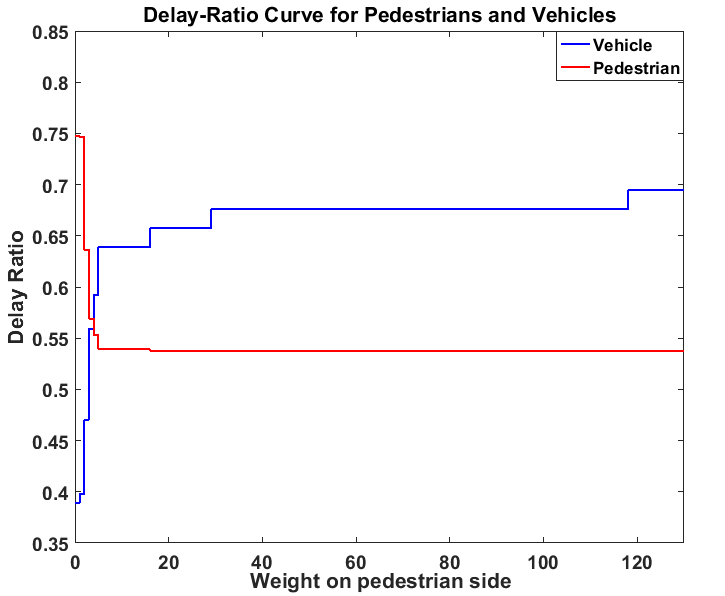}
	\caption{The change of delay-ratio curve under different weight of traffic composition from network 6$\times$6}
\label{fig:compare4}
\end{figure}

The regularity of switching frequency is investigated in Fig. \ref{fig:sf3}, which depicts the switching frequency of traffic signals at turning weight from network 3$\times$3 to network 10$\times$10. Additionally, the further impacts on the vehicle side are investigated from Fig. \ref{fig:compare1} to Fig. \ref{fig:compare4},
which illustrate the corresponding vehicle delay ratios and pedestrian delay ratios when continuously changing the weight on the pedestrian side for network 3$\times$3 to network 6$\times$6 with 30s prediction horizon respectively. The pedestrian delay ratio and the vehicle delay ratio are obtained by scaling the absolute delay time with the maximum delay ($P_{D}^{max}$ and $V_{D}^{max}$) from equation (\ref{eq:transPed}) and (\ref{eq:transVeh}) respectively.
In order to better summarize our primary findings, we give the succeeding texts below:
	\begin{enumerate}[label=\textbf{\arabic*},start=1]
		\item The states of vehicle traffic signals from the integrated model are exactly the same as the pure vehicle traffic model when the weight applying on the pedestrian side is 0.
		\item With the increase of the weight value, the vehicle signal states begin to change, as shown in Table \ref{integrate1} and \ref{integrate2} (the specific change signals are highlighted).
 Finally, the transform process is finished and all the vehicle signal states are exactly the same as the pedestrian traffic signal when optimize a pure pedestrian model, after that, no matter how heavy the weight is, the vehicle signal states will not change at all.
		\item The interesting phenomenon is that vehicle signals almost have a dramatic change at the initial several stages, then the change becomes more moderate. For example, in Table \ref{sf2},
25$\%$ and 40.63$\%$ of the traffic lights are highlighted at the second stage and the third stage shown in Table \ref{integrate2} when the weight is equal to 2 and 3, respectively, but after that, signals become increasingly prefers to maintain the current status unless more severe weights are applied. In order to adequately illustrate this phenomenon, three dimensional bar charts are depicted in Fig. \ref{fig:sf3}. The three coordinate axes denote network size, turning weight and associated switching frequency, respectively. Different network size is shown in different color. From the Fig. \ref{fig:sf3}, we can clearly see the trend that SF values are very higher at initial several stages, followed by a sharp decrease, then the values become  stable and small until switching to the individual pedestrian traffic light.
		\item The total vehicle delay presents a comprehensively rising trend and eventually becomes stabilized after completely converting to pure pedestrian signal, meanwhile the total pedestrian delay drifts steadily downward until the minimum pedestrian delay is reached. Also, all changes are in step form, which is shown
in Fig. \ref{fig:compare1} to Fig. \ref{fig:compare4}.
	\end{enumerate}

\begin{table}[!ht]
\centering
\caption{Pedestrian Impact On Vehicle side for Network 3 $\times$ 3 with 30s Prediction Horizon}
\label{integrate1}
\begin{tabular}{|c|c|c|c|c|c|c|c|c|c|c|c|c|}
\hline \hline
Weight m          & \multicolumn{2}{c|}{0}                                                             & \multicolumn{2}{c|}{2}                                                             & \multicolumn{2}{c|}{3}                                                             & \multicolumn{2}{c|}{4}                                                             & \multicolumn{2}{c|}{5}                                                             & \multicolumn{2}{c|}{8}                                                             \\ \hline
Hp=15s            & \multicolumn{2}{c|}{\begin{tabular}[c]{@{}c@{}}0 1 0\\ 1 0 1\\ 0 1 0\end{tabular}} & \multicolumn{2}{c|}{\begin{tabular}[c]{@{}c@{}}0 \textbf{0 1}\\ 1 \textbf{1 0}\\ \textbf{1 0 1}\end{tabular}} & \multicolumn{2}{c|}{\begin{tabular}[c]{@{}c@{}}0 0 1\\ 1 1 \textbf{1}\\ 1 0 1\end{tabular}} &
\multicolumn{2}{c|}{\begin{tabular}[c]{@{}c@{}}0 0 1\\ 1 1 1\\ 1 0 1\end{tabular}} &
\multicolumn{2}{c|}{\begin{tabular}[c]{@{}c@{}}0 0 1\\ 1 1 1\\ 1 0 1\end{tabular}} &
\multicolumn{2}{c|}{\begin{tabular}[c]{@{}c@{}}0 0 1\\ \textbf{0} 1 1\\ 1 0 1\end{tabular}} \\ \hline
Hp=30s            & \multicolumn{2}{c|}{\begin{tabular}[c]{@{}c@{}}1 0 1\\ 0 1 0\\ 1 0 1\end{tabular}} & \multicolumn{2}{c|}{\begin{tabular}[c]{@{}c@{}}1 0 \textbf{0}\\ 0 1 \textbf{1}\\ 1 \textbf{1 0}\end{tabular}} & \multicolumn{2}{c|}{\begin{tabular}[c]{@{}c@{}}\textbf{0} 0 0\\ 0 1 1\\ 1 1 \textbf{1}\end{tabular}} &
\multicolumn{2}{c|}{\begin{tabular}[c]{@{}c@{}}0 0 \textbf{1}\\ 0 1 1\\ 1 \textbf{0} 1\end{tabular}} &
\multicolumn{2}{c|}{\begin{tabular}[c]{@{}c@{}}0 0 1\\ \textbf{1} 1 1\\ 1 0 1\end{tabular}} &
\multicolumn{2}{c|}{\begin{tabular}[c]{@{}c@{}}0 0 1\\ \textbf{0} 1 1\\ 1 0 1\end{tabular}}
\\ \hline

\end{tabular}
\end{table}

\begin{table}[!ht]
\centering
\caption{Switching frequency of traffic signals at turning weight for Network 3 $\times$ 3 with 30s Prediction Horizon}
\renewcommand{\arraystretch}{1.4}
\label{sf1}
\resizebox{\columnwidth}{!}{%
\begin{tabular}[c]{|c|c|c|c|c|c|}
\hline
\hline
Turning Weight  & 2 & 3 & 4 & 5 & 8 \\ \hline
Switching  Frequency & 61.11\% & 16.67\% & 11.11\% & 5.56\% &11.11\% \\ \hline
\end{tabular}
}                                                                                                                        \\
\end{table}

\begin{table}[!ht]
\centering
\caption{Pedestrian Impact On Vehicle side for Network 4 $\times$ 4 with 30s Prediction Horizon}
\label{integrate2}
\begin{tabular}{|c|c|c|c|c|c|c|c|c|c|c|}
\hline \hline
\textbf{Weight m}                                            & \multicolumn{2}{c|}{0}                                                                             & \multicolumn{2}{c|}{2}                                                                             & \multicolumn{2}{c|}{3}                                                                             & \multicolumn{2}{c|}{5}                                                                             & \multicolumn{2}{c|}{6}                                                                             \\ \hline
Hp=15s
& \multicolumn{2}{c|}{\begin{tabular}[c]{@{}c@{}}0 1 0 1\\ 1 0 1 0\\ 0 1 0 1\\ 1 0 1 0\end{tabular}}
& \multicolumn{2}{c|}{\begin{tabular}[c]{@{}c@{}}\textbf{1} 1 0 1\\ 1 0 1 0\\ 0 \textbf{0} 0 1\\ \textbf{0} 0 1 0\end{tabular}}
& \multicolumn{2}{c|}{\begin{tabular}[c]{@{}c@{}}\textbf{0 0} 0 \textbf{0}\\ \textbf{0} 0 \textbf{0 1}\\ 0 0 0 1\\ 0 0 1 \textbf{1}\end{tabular}}
& \multicolumn{2}{c|}{\begin{tabular}[c]{@{}c@{}}0 0 0 0\\ 0 0 0 1\\ 0 0 0 1\\ 0 0 \textbf{0} 1\end{tabular}}
& \multicolumn{2}{c|}{\begin{tabular}[c]{@{}c@{}}0 0 0 0\\ 0 0 0 \textbf{0}\\ 0 0 0 1\\ 0 0 0 1\end{tabular}} \\ \hline
Hp=30s
& \multicolumn{2}{c|}{\begin{tabular}[c]{@{}c@{}}1 0 1 0\\ 0 1 0 1\\ 1 0 1 0\\ 0 1 0 1\end{tabular}}
& \multicolumn{2}{c|}{\begin{tabular}[c]{@{}c@{}}\textbf{0 1 0} 0\\ \textbf{1} 1 0 1\\ 1 0 1 0\\ 0 \textbf{0} 0 1\end{tabular}}
& \multicolumn{2}{c|}{\begin{tabular}[c]{@{}c@{}}0 \textbf{0} 0 0\\ \textbf{0 0} 0 \textbf{0}\\ \textbf{0} 0 1 \textbf{1}\\ 0 0 0 1\end{tabular}}
& \multicolumn{2}{c|}{\begin{tabular}[c]{@{}c@{}}0 0 0 0\\ 0 0 0 0\\ 0 0 1 1\\ 0 0 0 1\end{tabular}}
& \multicolumn{2}{c|}{\begin{tabular}[c]{@{}c@{}}0 0 0 0\\ 0 0 0 0\\ 0 0 1 1\\ 0 0 0 1\end{tabular}} \\ \hline                                                                                                                              \\ \hline \hline \hline
Weight m                                                     & \multicolumn{2}{c|}{11}                                                                            & \multicolumn{2}{c|}{14}                                                                            & \multicolumn{2}{c|}{30}                                                                            & \multicolumn{2}{c|}{}                                                                              & \multicolumn{2}{c|}{}                                                                              \\ \hline
Hp=15s
& \multicolumn{2}{c|}{\begin{tabular}[c]{@{}c@{}}0 0 0 0\\ 0 0 0 0\\ 0 0 0 1\\ 0 0 \textbf{1} 1\end{tabular}}
& \multicolumn{2}{c|}{\begin{tabular}[c]{@{}c@{}}0 \textbf{1} 0 0\\ 0 0 0 0\\ 0 0 0 1\\ 0 0 1 1\end{tabular}}
& \multicolumn{2}{c|}{\begin{tabular}[c]{@{}c@{}}0 1 0 0\\ 0 0 0 0\\ \textbf{1} 0 0 1\\ 0 0 1 1\end{tabular}}
& \multicolumn{2}{c|}{}
& \multicolumn{2}{c|}{}                                                                              \\ \hline
Hp=30s
& \multicolumn{2}{c|}{\begin{tabular}[c]{@{}c@{}}0 0 0 0\\ 0 0 0 0\\ 0 0 1 1\\ 0 0 \textbf{1} 1\end{tabular}}
& \multicolumn{2}{c|}{\begin{tabular}[c]{@{}c@{}}0 \textbf{1} 0 0\\ 0 0 0 0\\ 0 0 1 1\\ 0 0 1 1\end{tabular}}
& \multicolumn{2}{c|}{\begin{tabular}[c]{@{}c@{}}0 1 0 0\\ 0 0 0 0\\ \textbf{1} 0 1 1\\ 0 0 1 1\end{tabular}}
& \multicolumn{2}{c|}{}
& \multicolumn{2}{c|}{}
\\ \hline
\end{tabular}
\end{table}

\begin{table}[!ht]
\centering
\caption{Switching frequency of traffic signals at turning weight for Network 4 $\times$ 4 with 30s Prediction Horizon}
\renewcommand{\arraystretch}{1.4}
\label{sf2}
\resizebox{\columnwidth}{!}{%
\begin{tabular}[c]{|c|c|c|c|c|c|c|c|}
\hline
\hline
Turning Weight  & 2 & 3 & 5 & 6 & 11& 14 &30\\ \hline
Switch  Frequency & 25\% & 40.63\% & 3.13\% & 3.13\% & 6.25\% & 6.25\% & 6.25\%    \\ \hline
\end{tabular}
}                                                                                                                        \\
\end{table}

	Choosing a suitable weight is very difficult in order to balance both vehicles and pedestrians, especially when involving a large network and a long time horizon. Our findings give a way to quantitatively illustrate how pedestrians affect the vehicle network, which is expressed as a piecewise constant function. In our integrated model, the originally separated pedestrian network is now connected due to its coupling with the vehicle network, and a small change in the vehicle network could bring a significant impact on the whole system. In other words, the vehicle network has a much larger influence on the total vehicle network delay than the pedestrian network. For this reason, we need to put a larger weight value on the pedestrian side in order to match the influence of the vehicle network. For example, in Table \ref{integrate2}, the weight on the vehicle part is always 1 and the vehicle traffic signal begins to change until the weight on the pedestrian side is equal to 2, which means as long as the weight of the pedestrian side is between 0 and 2, its impact on the total traffic network delay cannot be obviously felt. The piecewise constant function of the vehicle network delay with respect to the weight values assigned to the pedestrian side shown from Fig. \ref{fig:compare1} to \ref{fig:compare4} clearly illustrates this fact.

\section{Conclusion} \label{sec:4}
In this paper we have proposed a pedestrian traffic light scheduling problem aiming to minimizing the total network delay and unhappiness. The pedestrian volume dynamic in each junction captures the characteristics of the conventional two-way crossing of the crosswalk topology, and the pedestrian flow dynamic is described by the source corner volume and the capacity of each crosswalk determined by the associated traffic signal setting and the hopping rule, which involve logical expression that can be converted into mixed integer linear constraints. By introducing two different evaluation criteria, the performance indexes are illustrated by pedestrian delay and pedestrian unhappiness, respectively. Then a vehicle traffic network is introduced and the integration of both the pedestrian delay and the vehicle delay is discussed. Our simulation results indicate that the pedestrian traffic light assignment can be handled efficiently when the total number of junctions and the prediction horizon are not too big. However, real-time scheduling cannot be guaranteed for a large urban traffic network, which can be effectively solved by DHS algorithm. Moreover, the switching frequency of the traffic signals show that unhappiness model do reduce the phenomenon, which frequently occur in delay model, that some pedestrians with few quantities but waiting for a long time. Additionally, our simulation also indicate that the influence of the pedestrian network on the total vehicle delay can be captured by a piecewise constant function, suggesting that the vehicle network may not always be sensitive to what happens in the pedestrian network, which could have some impact on the subsequent traffic light scheduling algorithm design for an integrated network, which shall be addressed in our future work. And current case study only consider vehicle network with one-way vehicle flows, more complex phases involving turning flows should also be taken into account for integrated model in our future work.
\bibliography{ref}

\begin{thebibliography}{10}

\bibitem{allsop1992evolving}
RE~Allsop.
\newblock Evolving application of mathematical optimisation in design and
  operation of individual signal-controlled road junctions.
\newblock In {\em Mathematics in Transport Planning and Control (Institute of
  Mathematics \& ITS Applications Conference Series 38)}, 1992.

\bibitem{allsop1971sigset}
Richard~E Allsop.
\newblock Sigset: a computer program for calculating traffic signal settings.
\newblock {\em Traffic Engineering \& Control}, 1971.

\bibitem{allsop1976sigcap}
Richard~E Allsop.
\newblock Sigcap: A computer program for assessing the traffic capacity of
  signal-controlled road junctions.
\newblock {\em Traffic Engineering \& Control}, 17(Analytic), 1976.

\bibitem{boillot1992optimal}
F~Boillot.
\newblock Optimal signal control of urban traffic networks.
\newblock In {\em International Conference on Road Traffic Monitoring and
  Control (6th: 1992: London, England). 6th International Conference on Road
  Traffic Monitoring and Control}, 1992.

\bibitem{boudet2009pedestrian}
Laurence Boudet and Sophie Midenet.
\newblock Pedestrian crossing detection based on evidential fusion of
  video-sensors.
\newblock {\em Transportation research part C: emerging technologies},
  17(5):484--497, 2009.

\bibitem{daganzo1994cell}
Carlos~F Daganzo.
\newblock The cell transmission model: A dynamic representation of highway
  traffic consistent with the hydrodynamic theory.
\newblock {\em Transportation Research Part B: Methodological}, 28(4):269--287,
  1994.

\bibitem{gaisbauer2008wayfinding}
Christian Gaisbauer and Andrew~U Frank.
\newblock Wayfinding model for pedestrian navigation.
\newblock In {\em AGILE 2008 Conference-Taking Geo-information Science One Step
  Further, University of Girona, Spain}, 2008.

\bibitem{gallivan1988optimising}
Stephen Gallivan and Benjamin Heydecker.
\newblock Optimising the control performance of traffic signals at a single
  junction.
\newblock {\em Transportation Research Part B: Methodological}, 22(5):357--370,
  1988.

\bibitem{gaarder1989pedestrian}
Per G{\aa}rder.
\newblock Pedestrian safety at traffic signals: a study carried out with the
  help of a traffic conflicts technique.
\newblock {\em Accident Analysis \& Prevention}, 21(5):435--444, 1989.

\bibitem{gartner1981versatile}
Nathan~H Gartner.
\newblock A versatile program for setting signals on arteries and triangular
  networks john dc little* mark d. kelson.
\newblock 1981.

\bibitem{gartner1983opac}
Nathan~H Gartner.
\newblock {\em OPAC: A demand-responsive strategy for traffic signal control}.
\newblock Number 906. 1983.

\bibitem{gazis1963oversaturated}
Denos~C Gazis and Renfrey~Burnard Potts.
\newblock The oversaturated intersection.
\newblock Technical report, 1963.

\bibitem{geem2001new}
Zong~Woo Geem, Joong~Hoon Kim, and GV~Loganathan.
\newblock A new heuristic optimization algorithm: harmony search.
\newblock {\em Simulation}, 76(2):60--68, 2001.

\bibitem{optimization2012gurobi}
Gurobi, Optimization, et~al.
\newblock Gurobi optimizer reference manual.
\newblock {\em URL: http://www. gurobi. com}, 2:1--3, 2012.

\bibitem{han2014continuum}
Ke~Han, Vikash~V Gayah, Benedetto Piccoli, Terry~L Friesz, and Tao Yao.
\newblock On the continuum approximation of the on-and-off signal control on
  dynamic traffic networks.
\newblock {\em Transportation Research Part B: Methodological}, 61:73--97,
  2014.

\bibitem{henry1984prodyn}
Jean-Jacques Henry, Jean~Loup Farges, and J~Tuffal.
\newblock The prodyn real time traffic algorithm.
\newblock In {\em IFACIFIPIFORS conference on control in}, 1984.

\bibitem{heydecker1987calculation}
Benjamin~G Heydecker and Ian~W Dudgeon.
\newblock Calculation of signal settings to minimise delay at a junction.
\newblock {\em Transportation and Traffic Theory}, 1987.

\bibitem{hunt1982scoot}
PB~Hunt, DI~Robertson, RD~Bretherton, and M~Cr Royle.
\newblock The scoot on-line traffic signal optimisation technique.
\newblock {\em Traffic Engineering \& Control}, 23(4), 1982.

\bibitem{improta1984control}
G~Improta and GE~Cantarella.
\newblock Control system design for an individual signalized junction.
\newblock {\em Transportation Research Part B: Methodological}, 18(2):147--167,
  1984.

\bibitem{ishaque2005multimodal}
Muhammad Ishaque and Robert Noland.
\newblock Multimodal microsimulation of vehicle and pedestrian signal timings.
\newblock {\em Transportation Research Record: Journal of the Transportation
  Research Board}, (1939):107--114, 2005.

\bibitem{ishaque2007trade}
Muhammad~Moazzam Ishaque and Robert~B Noland.
\newblock Trade-offs between vehicular and pedestrian traffic using
  micro-simulation methods.
\newblock {\em Transport Policy}, 14(2):124--138, 2007.

\bibitem{kasemsuppakorn2013pedestrian}
Piyawan Kasemsuppakorn and Hassan~A Karimi.
\newblock A pedestrian network construction algorithm based on multiple gps
  traces.
\newblock {\em Transportation research part C: emerging technologies},
  26:285--300, 2013.

\bibitem{lam1997integrated}
William~HK Lam, Antonio~CK Poon, and Gregory~KS Mung.
\newblock Integrated model for lane-use and signal-phase designs.
\newblock {\em Journal of transportation engineering}, 123(2):114--122, 1997.

\bibitem{lin2001enhanced}
Wei-Hua Lin.
\newblock An enhanced 0-1 mixed integer lp formulation for the traffic signal
  problem.
\newblock In {\em Intelligent Transportation Systems, 2001. Proceedings. 2001
  IEEE}, pages 189--194. IEEE, 2001.

\bibitem{little1966synchronization}
John~DC Little.
\newblock The synchronization of traffic signals by mixed-integer linear
  programming.
\newblock {\em Operations Research}, 14(4):568--594, 1966.

\bibitem{lo1999novel}
Hong~K Lo.
\newblock A novel traffic signal control formulation.
\newblock {\em Transportation Research Part A: Policy and Practice},
  33(6):433--448, 1999.

\bibitem{ma2015optimization}
Wanjing Ma, Dabin Liao, Yue Liu, and Hong~Kam Lo.
\newblock Optimization of pedestrian phase patterns and signal timings for
  isolated intersection.
\newblock {\em Transportation Research Part C: Emerging Technologies},
  58:502--514, 2015.

\bibitem{ma2014optimization}
Wanjing Ma, Yue Liu, and K~Larry Head.
\newblock Optimization of pedestrian phase patterns at signalized
  intersections: a multi-objective approach.
\newblock {\em Journal of advanced transportation}, 48(8):1138--1152, 2014.

\bibitem{manual2010hcm2010}
Highway~Capacity Manual.
\newblock Hcm2010.
\newblock {\em Transportation Research Board, National Research Council,
  Washington, DC}, 2010.

\bibitem{miller1963computer}
Alan~J Miller.
\newblock A computer control system for traffic networks.
\newblock 1963.

\bibitem{perkins1968traffic}
Stuart~R Perkins and Joseph~L Harris.
\newblock Traffic conflict characteristics-accident potential at intersections.
\newblock {\em Highway Research Record}, (225), 1968.

\bibitem{robertson1900tansyt}
Dennis~I Robertson.
\newblock 'tansyt'method for area traffic control.
\newblock {\em Traffic Engineering \& Control}, 8(8), 1900.

\bibitem{sen1997controlled}
Suvrajeet Sen and K~Larry Head.
\newblock Controlled optimization of phases at an intersection.
\newblock {\em Transportation science}, 31(1):5--17, 1997.

\bibitem{tarko1995accident}
Andrzej Tarko and Marian Tracz.
\newblock Accident prediction models for signalized crosswalks.
\newblock {\em Safety science}, 19(2):109--118, 1995.

\bibitem{virkler1998scramble}
Mark Virkler.
\newblock Scramble and crosswalk signal timing.
\newblock {\em Transportation Research Record: Journal of the Transportation
  Research Board}, (1636):83--87, 1998.

\bibitem{webster1958traffic}
Fo~Vo Webster.
\newblock Traffic signal settings.
\newblock Technical report, 1958.

\bibitem{zhang2015urban}
Yicheng Zhang, Rong Su, and Kaizhou Gao.
\newblock Urban road traffic light real-time scheduling.
\newblock In {\em 2015 54th IEEE Conference on Decision and Control (CDC)},
  pages 2810--2815. IEEE, 2015.

\end{thebibliography}
\bibliographystyle{ieeetran}
\end{document}